\documentclass{IEEEtran4}
\ifCLASSINFOpdf
   \usepackage[pdftex]{graphicx}
\else
   \usepackage[dvips]{graphicx}
\fi
\usepackage{amsmath,amsfonts,amssymb,amsthm}
\usepackage{bbm}
\usepackage{algorithm}
\usepackage{algpseudocode}
\usepackage{booktabs}

\usepackage[dvipsnames]{xcolor}
\usepackage{todonotes}

\newcommand{\mt}[1]{\todo[inline,color=blue!30]{[MT] #1}}

\newcommand{\raproxy}{{\sc RA-E2ELR}}

\ifCLASSOPTIONcompsoc
 \usepackage[caption=false,font=normalsize,labelfont=sf,textfont=sf]{subfig}
\else
 \usepackage[caption=false,font=footnotesize]{subfig}
\fi
\usepackage{float}
\floatstyle{ruled}
\newfloat{model}{thp}{lop}
\floatname{model}{Model}

\newcommand{\ubar}[1]{\underline{#1}}

\newcommand{\e}{\mathbf{e}}  % vector of all ones

\newcommand{\x}{\mathbf{x}}

\newcommand{\Mpb}{M_{\text{pb}}}

\newcommand{\Mth}{M_{\text{th}}}

\newcommand{\pd}{\mathbf{d}}
\newcommand{\pfmax}{\mathbf{\bar{f}}}
\newcommand{\pfmin}{\mathbf{\ubar{f}}}
\newcommand{\pg}{\mathbf{p}}
\newcommand{\pgmin}{\mathbf{\ubar{p}}}
\newcommand{\pgmax}{\mathbf{\bar{p}}}

\newcommand{\xith}{\xi_{\text{th}}}

\newcommand{\btheta}{\boldsymbol{\Theta}}

\newcommand{\etaup}{\eta^\uparrow}
\newcommand{\etadn}{\eta^\downarrow}

\newcommand{\rup}{R^{\uparrow}}
\newcommand{\rdn}{R^{\downarrow}}

\makeatletter
\let\old@ps@headings\ps@headings
\let\old@ps@IEEEtitlepagestyle\ps@IEEEtitlepagestyle
\def\psccfooter#1{%
    \def\ps@headings{%
        \old@ps@headings%
        \def\@oddfoot{\strut\hfill#1\hfill\strut}%
        \def\@evenfoot{\strut\hfill#1\hfill\strut}%
    }%
    \def\ps@IEEEtitlepagestyle{%
        \old@ps@IEEEtitlepagestyle%
        \def\@oddfoot{\strut\hfill#1\hfill\strut}%
        \def\@evenfoot{\strut\hfill#1\hfill\strut}%
    }%
    \ps@headings%
}
\makeatother

\DeclareMathOperator*{\argmin}{arg-min}

\begin{document}
\title{Real-Time Risk Analysis with Optimization Proxies}
% % The code below cheats with margins to have a one-line title (but it overflows the margins)
% \title{\mbox{\hspace{-0.5em}Real-Time Risk Assessment with Optimization Proxies}}

%% To specify the authors when (number of affiliations <= 2)
\author{
\IEEEauthorblockN{Wenbo Chen, Mathieu Tanneau, Pascal Van Hentenryck}
\IEEEauthorblockA{Georgia Institute of Technology \\
wenbo.chen@gatech.edu, \{mathieu.tanneau, pascal.vanhentenryck\}@isye.gatech.edu}
}

%% To specify the authors when (number of affiliations > 2)
% \author{\IEEEauthorblockN{Author n.1\IEEEauthorrefmark{1},
% Author n.2\IEEEauthorrefmark{2},
% Author n.3\IEEEauthorrefmark{3}, 
% Author n.4\IEEEauthorrefmark{3} and
% Author n.5\IEEEauthorrefmark{4}}
% \IEEEauthorblockA{\IEEEauthorrefmark{1} Department Name of Organization A\\
% Name of the organization A,
% Address A\\ Emails if wanted}
% \IEEEauthorblockA{\IEEEauthorrefmark{2} Department Name of Organization B\\
% Name of the organization B,
% Address B\\ Emails if wanted}
% \IEEEauthorblockA{\IEEEauthorrefmark{3} Department Name of Organization C\\
% Name of the organization C,
% Address C\\ Emails if wanted}
% \IEEEauthorblockA{\IEEEauthorrefmark{4}Department Name of Organization D\\
% Name of the organization D,
% Address D\\ Emails if wanted}
% }

\maketitle

\begin{abstract}
    The increasing penetration of renewable generation and distributed energy resources requires new operating practices for power systems, wherein risk is explicitly quantified and managed.
    However, traditional risk-assessment frameworks are not fast enough for real-time operations, because they require numerous simulations, each of which requires solving multiple economic dispatch problems sequentially.
    The paper addresses this computational challenge by proposing proxy-based risk assessment, wherein optimization proxies are trained to learn the input-to-output mapping of an economic dispatch optimization solver.
    Once trained, the proxies make predictions in milliseconds, thereby enabling real-time risk assessment.
    The paper leverages self-supervised learning and end-to-end-feasible architecture to achieve high-quality sequential predictions.
    Numerical experiments on large systems demonstrate the scalability and accuracy of the proposed approach.
    % Power systems are experiencing growing operational uncertainty caused by the rise of renewable generators and distributed energy resources.
    % Therefore, operators have to continuously monitor risk in real-time by quickly assessing the system's behavior under changing load and renewable conditions.
    % However, solving an optimization problem like Security-Constrained Economic Dispatch (SCED) for each such scenario is computationally intractable given the real-time constraints. 
    % To address this challenge, this paper proposes a learning-based risk assessment pipeline that replaces the optimization solver with machine learning surrogates, i.e., optimization proxies.
    % These proxies can generate dispatch decisions that are both feasible and close to optimal within milliseconds, significantly reducing the time required for risk simulation from 15 minutes to the order of milliseconds.
    % Numerical experiments conducted on the Midcontinent Independent System Operator (MISO) system validate the effectiveness of the proposed ML-based risk assessment approach.
    % It successfully provides fine-grained risk measurements for large-scale power grids in industrial settings.
\end{abstract}

\begin{IEEEkeywords}
    Risk assessment, optimization proxies
\end{IEEEkeywords}

\thanksto{\noindent This research was partly supported by NSF award 2112533 and ARPA-E PERFORM award AR0001136.}

\section{Introduction}
\label{sec:introduction}

The growing penetration of distributed energy resources (DERs) and of
renewable generation, especially wind and solar generation, is causing
a fundamental change of paradigm for power systems operations.
Traditional approaches, where uncertainty is low and can be managed
through reserve products, are no longer viable, in presence of
significant wind and solar generation, and less predictable demand
patterns due to DER adoption. Modern power grids thus require new
operating practices that explicitly assess operational risk in real
time.

Traditional approaches to risk assessment require running large
numbers (typically in the thousands) of Monte-Carlo (MC) simulations.
Risk and reliability metrics are then computed based on the outputs of
these simulations \cite{stover2023reliability}.  In power systems,
each MC simulation evaluates the behavior of the grid, e.g., generator
setpoints and power flows, under a possible realization of load and
renewable output, obtained from solving multiple economic dispatch
problems sequentially (one per time step).  Nevertheless, the
computational cost of optimization, combined with the large number of
MC simulations needed to obtain an accurate risk estimate, renders
this workflow intractable for real-time risk assessment.

\begin{figure}[!t]
        \centering
        \subfloat[Optimization-based Risk Assessment.]{
           \label{fig:risk_opt}
           \includegraphics[width=0.95\columnwidth]{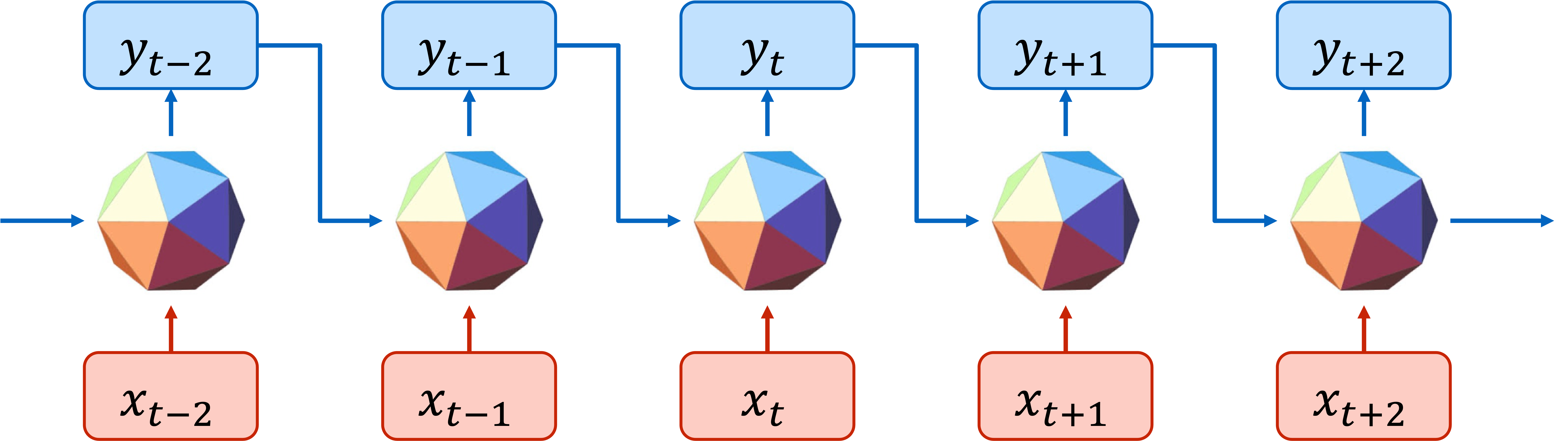}
        }\\[1em]
        \subfloat[Risk Assessment with E2ELR Proxies.]{
           \label{fig:risk_ml}
           \includegraphics[width=0.95\columnwidth]{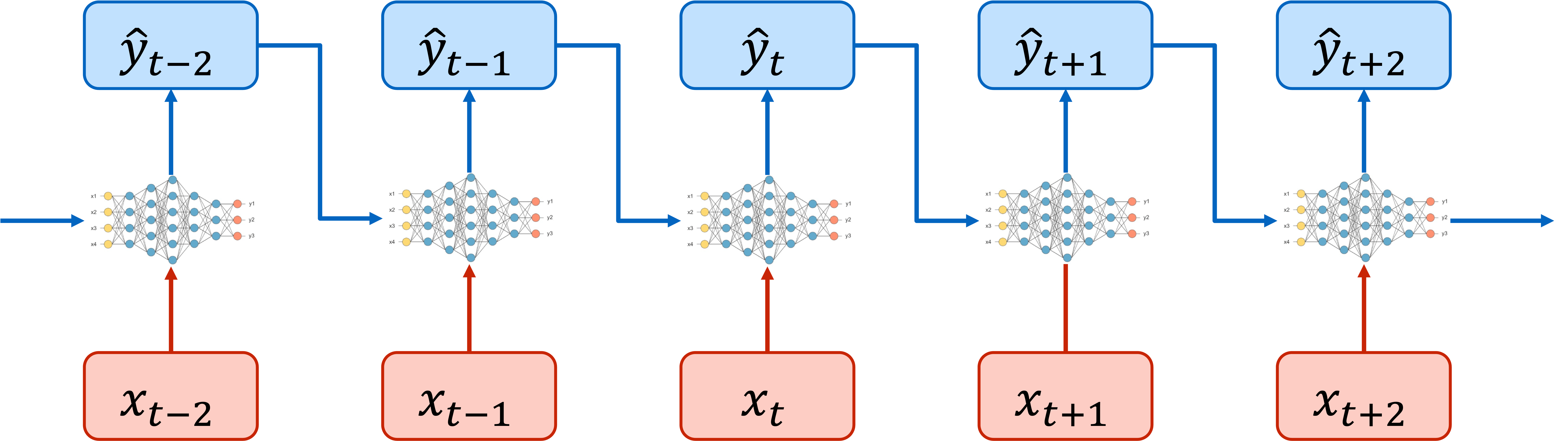}
        }
        \caption{Illustration of the optimization-based (top) and \raproxy{} MC simulation pipeline. At each step $t$, generator setpoints $y_{t}$ are decided based on current system conditions $x_{t}$ and previous setpoints $y_{t-1}$.}
        % \caption{a) Run risk simulation by solving multiple optimization models sequentially. It takes more than 15 mins to run on an industrial-size system like MISO. b) Run risk simulation with optimization proxies sequentially. It takes less than 5 milliseconds to run the whole simulation.}
        \label{fig:risk_model}
        \vspace*{-1em}
    \end{figure}

The paper proposes a new methodology, Risk Assessment with End-to-End
Learning and Repair (\raproxy{}), to address this computational
challenge.  \raproxy{} replaces costly optimization solvers with fast
E2ELR optimization proxies, i.e., differentiable programs that are
guaranteed to produce near-optimal feasible solutions to economic
dispatch problems \cite{chen2023end}. \raproxy{} is sketched in Figure
\ref{fig:risk_model} and is capable of assessing risk in electricity
markets whose real-time operations are organized around economic
dispatch processes that co-optimize energy and reserves. The paper
reports numerical experiments on large-scale power grids with
thousands of buses, demonstrating the scalability and fidelity of
\raproxy{}. To the best of the authors' knowledge, \raproxy{} is the
first methodology to real-time risk assessment for the US
market-clearing pipeline.
    
    % \textcolor{blue}{[MT: done until here]}

    % Power systems are expected to incorporate substantial amounts of renewable generation and distributed energy resources in the near future \cite{Francis_2021}.
    % This brings increasing stochasticity in load and renewable generation, resulting in more significant prediction errors.
    % In current practice, system operators handle uncertainty in daily operations by acquiring excess unused capacity i.e., reserve products.
    % However, these implicit risk management approaches may no longer be adequate in operational regimes characterized by higher uncertainty and narrower safety margins.
    % As a result, explicitly modeling the operational risk is critical to inform system operators to enable proactive risk management \cite{stover2023reliability}.

    % However, in the context of real-time system operation, computing the explicit risk assessment needs solving a large number Security-Constrained Economic Dispatch (SCED) instances corresponding to different scenarios including loads, wind, and solar generation. Figure~\ref{fig:risk_opt} illustrates one trajectory of the risk assessment, which involves solving $T$ SCED instances sequentially. 
    % It already takes more than 15 mins to run on an industrial-size system like MISO.
    % Computing the risk metrics such as conditional expectations and probabilities of failure needs to collect thousands of trajectories.
    % Thus it is computationally intractable for real-time system operation.

\section{Literature Review}

\subsection{Risk Assessment in Power Systems}

A number of reliability metrics have been proposed to quantify the
risk in power systems explicitly such as Lack of ramp probability
(LORP) \cite{thatte2015metric}, insufficient ramping resource
expectation (IRRE) \cite{lannoye2012evaluation}, loss of load
expectation (LOLE) and expected unserved energy (EUE) \cite{lole}.
More recently, Stover \& al \cite{stover2023reliability} proposed
three levels of risk metrics for a quantity of interest: conditional
value at risk, reliability, and risk.  Computing those risk metrics is
based on Monte Carlo (MC) sampling, which requires solving
optimization instances for thousands of scenarios.  These MC
simulations are computationally expensive and intractable in real-time
settings. This paper shows how to compute, in real-time, high-fidelity
approximations of these risk metrics by leveraging E2ELR optimization
proxies.

    % \begin{itemize}
    %     \item \cite{stover2023reliability} proposes comprehensive risk metrics at three levels: conditional expectation, reliability, and risk.
    %     \item \cite{stover2022just} ML-based risk simulation for power system at system level.
    %     \item  \cite{bo2023dnn} DNN-based risk simulation for multi-state series-parallel systems considering semi-Markov process.
    % \end{itemize}

\subsection{Optimization Proxies for Sequential OPF}

Some recent papers
\cite{zhou2020data,yan2020real,sayed2022feasibility,sayed2023optimal,zeng2022physics}
focus on developing optimization proxies for multi-period ACOPF
problems. They model the proxies with deep neural networks, typically
combined with repair steps for feasibility.  The proxies are trained
under the framework of Reinforcement Learning (RL) with different
algorithms.  For instance, in \cite{zhou2020data}, proxies are
pretrained with imitation learning and then improved using proximal
policy optimization.  The results are reported for systems with up to
200 buses.  In \cite{yan2020real}, the authors train proxies using
deep deterministic policy gradient.  The reward function is augmented
with constraint violations to improve the feasibility.  The results
are reported on the 118-bus system. {\em Those methodologies cannot
  guarantee feasibility and thus cannot be used for risk assessments}.
    
Subsequent results addressed feasibility issues by designing complex
optimization-based feasibility layers such as power flow solvers and
convex-concave procedures \cite{sayed2022feasibility}, or holomorphic
embedding \cite{sayed2023optimal}; they also use more advanced RL
algorithms such as the soft actor-critic framework.  Because of the
complexity of their optimization-based repair steps, {\em the
  improvements due to the proxies compared to the optimization
  algorithms are not sufficient for real-time risk assessment}.  The
results are only reported up to systems with 300 buses.

    % Different from those works, the proxies of this work are built on the End-to-End Learning and Repair (E2ELR) architecture in \cite{chen2023end}, where the model could be trained using self-supervised learning, rather than sample-expensive reinforcement learning. The feasibility is also guaranteed by designing efficient and close-form feasibility layers, rather than computational-expensive optimization-based projections.
    % \cite{zeng2022physics} integrate deep reinforcement learning and physics-constrained adversary example generation for multi-period AC-SCOPF.

\section{Operational Risk Assessment in Power Grids} 
\label{sec:formulation}

Consider a power grid with buses $\mathcal{N} = \{1, ..., N\}$, and
branches $\mathcal{E} \, {=} \, \{1, ..., E\}$.  Assume, without loss
of generality, that the costs are linear and that exactly one
generator and one load are attached to each bus.

\subsection{The Economic Dispatch Formulation}
\label{sec:formulation:ED}   

    \begin{model}[!t]
        \caption{The Economic Dispatch Model.}
        \label{model:ED}
        \textbf{Input}: Current load $\pd_{t}$, previous generation dispatch $\pg_{t-1}$\\
        \textbf{Output}: Generation dispatch $\pg_{t}$
        \begin{subequations}
        \label{eq:ED}
        \begin{align}
            \min_{\pg_{t}} \quad 
            & \label{eq:ED:objective}
                c^{\top}\pg_{t} + \Mth \| \xith \|_{1} \\
            \text{s.t.} \quad
            & \label{eq:ED:power_balance}
               \e^{\top}\pg_{t} = \e^{\top}\pd_{t}\\
            & \label{eq:ED:min_max_limits}
               \pgmin \leq \pg_{t} \leq \pgmax\\
            % & \label{eq:ED:ramp_dn}
            %    \pg_{t} \geq \pg_{t-1} - \rdn\\
            % & \label{eq:ED:ramp_up}
            %    \pg_{t} \leq \pg_{t-1} + \rup\\
            & \label{eq:ED:ramp}
               \pg_{t-1} - \rdn \leq \pg_{t} \leq \pg_{t-1} + \rup\\
            & \label{eq:ED:thermal_constraints}
               \pfmin - \xith \leq \Phi (\pg_{t} - \pd_{t}) \leq \pfmax + \xith
        \end{align}
        \end{subequations}
    \end{model}    
    
    Model \ref{model:ED} presents the Economic Dispatch (ED) model to
    be solved at time step $t \, {\in} \, \mathbb{Z}$. For simplicity,
    the formulation omits the reserve constraints but \raproxy{} has
    no difficulty in handling them \cite{chen2023end}.  The duration
    of each time step is assumed to be $1$; in the US, ED is solved
    every 5 minutes to clear the real-time electricity market.  The ED
    problem at time $t$ takes, as inputs, the current demand $\pd_{t}$
    and the prior generation dispatch $\pg_{t-1}$, and outputs the
    generation dispatch $\pg_{t}$.  For generators, the minimum
    output, maximum output, upward and downward ramping capacity are
    denoted by $\pgmin$, $\pgmax$, $\rup$ and $\rdn$, respectively.
    The Power Transfer Distribution Factor (PTDF) matrix is denoted by
    $\Phi \, {\in} \, \mathbb{R}^{E \times N}$, and the lower and
    upper thermal limits are denoted by $\pfmin$ and $\pfmax$.
    
    The ED objective \eqref{eq:ED:objective} minimizes the sum of
    total generation costs and thermal violation penalties.
    Constraint \eqref{eq:ED:power_balance} enforces power balance in
    the system.  Constraints \eqref{eq:ED:min_max_limits} enforce
    generation minimum and maximum limits, and constraints
    \eqref{eq:ED:ramp} enforce ramping constraints.  Constraints
    \eqref{eq:ED:min_max_limits}--\eqref{eq:ED:ramp} can be
    combined into
    \begin{align}
        \label{eq:ED:bounds}
        \pgmin_{t} \, {\leq} \, \pg_{t} \, \leq \pgmax_{t},
    \end{align}
    where $\pgmin_{t} {=} \, \max \left( \pgmin, \pg_{t-1} {-} R^{\downarrow} \right)$ and $\pgmax_{t} {=} \, \min \left( \pgmax, \pg_{t-1} {+} R^{\uparrow} \right)$.
    Finally, constraints \eqref{eq:ED:thermal_constraints} express thermal constraints using a PTDF formulation.
    Following the practice of US system operators \cite{MISO_BPM_002}, thermal constraints are treated as soft, i.e., thermal violations $\xith$ are penalized using a high price $\Mth$.

\subsection{The Operational Risk Assessment Framework}
\label{sec:formulation:risk}
    
    Principled risk assessment is based on executing large numbers of
    MC simulations \cite{stover2023reliability}.  Each individual MC
    simulations emulates the system behavior under given (net) load
    conditions.  The results of multiple MC simulations are then
    aggregated into risk and reliability metrics.

    Formally, consider a simulation horizon $\mathcal{T} \, {=} \,
    \{1, ..., T\}$.  In a day-ahead setting, $\mathcal{T}$ may
    correspond to each hour of the operating day; in a real-time
    setting, $\mathcal{T}$ may capture the subsequent few hours at a
    5-minute granularity.  A set of scenarios $\mathcal{S} \, {=} \,
    \{1, ..., S\}$, typically generated by a probabilistic forecasting
    model, is available.  Algorithm \ref{alg:risk} formalizes the
    simulation algorithm, which is also illustrated in Figure
    \ref{fig:risk_opt}.  Each simulation $s$ takes as input a sequence
    of load $(\pd^{s}_{1}, ..., \pd^{s}_{T})$, where $\pd^{s}_{t}$
    denotes the vector of loads a time $t$ and for scenario $s$.  At
    each time $t$, the minimum and maximum limits $\pgmin^{s}_{t},
    \pgmax^{s}_{t}$ are updated.  Then, current generation dispatches
    $\pg^{s}_{t}$ are computed by solving the ED problem in Model
    \ref{model:ED}.  Finally, the simulation returns the sequence of
    generation dispatches $(\pg^{s}_{1}, ..., \pg^{s}_{T})$.

    \begin{algorithm}[!t]
        \caption{Monte-Carlo Simulation}
        \label{alg:risk}
        \begin{algorithmic}[1]
        \Require Demand scenario $(\pd^{s}_{1}, ..., \pd^{s}_{T})$
            \For{$t = 1 \cdots T$}
                \State $\pgmax^{s}_t, \pgmin^{s}_t \gets \max \left( \pgmin, \pg^{s}_{t-1} {-} R^{\downarrow} \right), \min \left( \pgmax, \pg^{s}_{t-1} {+} R^{\uparrow} \right)$
                \State $\pg^{s}_t \gets \text{ED}(\pd^{s}_t, \pgmax^{s}_t, \pgmin^{s}_t)$
            \EndFor{}
            \State \textbf{Return} $(\pg^{s}_{1}, ..., \pg^{s}_{T})$
        \end{algorithmic}
    \end{algorithm}

    Risk metrics are computed by aggregating the outputs of multiple MC simulations, typically hundreds or thousands.
    This work considers the three levels of risk metrics proposed in \cite{stover2023reliability}: conditional value at risk, probability of failure, and risk.
    Let $\mathcal{Q}$ be a so-called \emph{quantity of interest} (QoI), whose value depends on the state of the system.
    Namely, denote by 
    \begin{align}
        \label{eq:risk:QOI}
        Q^{s}_{t} = \mathcal{Q}(\pd^{s}_{t}, \pg^{s}_{t})
    \end{align}
    the QoI value at time $t$ in scenario $s$.  For instance, if
    $\mathcal{Q}$ is the system power imbalance, then $Q^{s}_{t} =
    \e^{\top}(\pg^{s}_{t} - \pd^{s}_{t})$.
    
    \subsubsection{Conditional Value at Risk (CVAR)}
        Let $\alpha \, {\in} \, [0, 1]$ be a given significance level, typically chosen by a system operator. 
        The CVAR metric for time $t$ is
        \begin{align}
            \label{eq:risk_metrics:condition_expectation}
            \text{CE}(\mathcal{Q}_{t}) 
                = \mathbb{E}[\mathcal{Q}_{t} | \mathcal{Q}_{t} \geq \mathcal{Q}_{t}(\alpha)] 
                \approx
                \frac{
                    \displaystyle \sum_{s \in \mathcal{S}} Q^{s}_{t} \cdot \mathbbm{1}_{Q^{s}_{t} \geq Q_{t}(\alpha)}
                }{
                    \displaystyle \sum_{s \in \mathcal{S}} \mathbbm{1}_{Q^{s}_{t} \geq Q_{t}(\alpha)}
                },
        \end{align}
        where $\mathbbm{1}$ denotes the indicator function and $Q_{t}(\alpha)$ denote the $\alpha$ quantile of $(Q^{1}_{t}, ..., Q^{S}_{t})$.
        Note that the CVAR metric can be adapted to capture the left tail of the distribution \cite{stover2023reliability}.
    
    \subsubsection{Probability of Failure}
        To assess system reliability, the second metric considers the probability of failure.
        Namely, let $\bar{Q}$ be a desired reliability threshold, i.e., the system is considered to fail if $\mathcal{Q} \geq \bar{Q}$.
        Then, the probability of failure at time $t$ is given by
        \begin{align}
            \label{eq:risk_metrics:reliability}
            \mathbb{P}(\mathcal{Q}_{t} \geq \bar{Q}) \approx \frac{1}{S} \sum_{s \in \mathcal{S}} \mathbbm{1}_{Q^{s}_{t} \geq \bar{Q}}.
        \end{align}
        Note that reliability can be measured for the entire system, or for individual components.
        % The reliability considers the probability of occurrence of these events related to the quantity of interest.
        % Given a quantity of interest $R$, a set of $N$ samples $\{R_n\}_{n=1}^N$ and $\bar{R}$ denotes the desired reliability threshold, the probability of the failure event i.e., $R < \bar{R}$ is defined as:
    
    \subsubsection{Risk}
        Risk quantifies the monetary consequences of failure.
        The risk at time $t$ is given by
        \begin{align}
            \label{eq:risk_metrics:risk}
            \text{Risk}(\mathcal{Q}_{t}) = \mathbb{E}[c(\mathcal{Q}_{t})] \approx \frac{1}{S}\sum_{s \in \mathcal{S}} c(Q^{s}_{t}),
        \end{align}
        where $c$ denotes a cost function.
        For instance, if the QoI $\mathcal{Q}$ denotes system power imbalance, a natural choice for $c$ could be $c \, {=} \, \text{VOLL} {\times} \mathcal{Q}$, where VOLL is the Value of Lost Load.
        Compared to the previous two metrics, risk captures the magnitude of adverse events and their financial consequences.
        
        % Given a quantity of interest $R$, a set of $N$ samples $\{R_n\}_{n=1}^N$ and $c(.)$ denotes the cost function of the event failure, the risk of the event failure is defined as:
        % where $c(\cdot)$ is typically defined as a piece-wise linear function in power systems.

\section{Risk Assessment With \raproxy{}}
\label{sec:methodology}

\newcommand{\z}{\mathbf{z}}
\newcommand{\pghat}{\mathbf{\hat{p}}}
\newcommand{\pgfeas}{\mathbf{\tilde{p}}}

The main limitation of existing risk assessment frameworks is their
computational cost.  Indeed, running thousands of simulations requires
solving tens to hundreds of thousands of ED problems.  This approach
is not tractable, especially as more resources are integrated into the
grid and the complexity of ED problems increases in the future.

\subsection{Proxy-Based Risk Assessment}
\label{sc:ML:proxy_risk}

To enable fast risk assessment, the paper replaces costly optimization
solvers with fast optimization proxies, i.e., {\em differentiable
  programs that approximate the input-output mapping of optimization
  solvers}.  Figure \ref{fig:risk_ml} illustrates the proposed risk
assessment framework, which replaces the ED optimization in step 3 of
Algorithm \ref{alg:risk} with an optimization proxy.  At each
time step $t$, the ED proxy takes as input the ED problem data
$\mathbf{x}_{t}$, and outputs a near-optimal solution
$\mathbf{\hat{y}}_{t}$ in milliseconds.  It is important to note that,
in the present setting, the proxy prediction at time $t$ is part of
the input of the subsequent prediction.  The sequential nature of the
prediction task is in stark contrast with existing literature on
optimization proxies, which consider individual ED instances.

The proxy-based risk assessment task can therefore be stated as
follows.  Given scenario $s$ with load $(\pd^{s}_{1}, ...,
\pd^{s}_{T})$, the model should output a sequence $(\pghat^{s}_{1},
..., \pghat^{s}_{T})$ of generator dispatches that is as close as
possible to the sequence $(\pg^{s}_{1}, ..., \pg^{s}_{T})$ produced by
an optimization solver.
    
Risk metrics can then be evaluated based on the predicted sequence
$(\pghat^{s}_{1}, ..., \pghat^{s}_{T})$.  Given a QoI $\mathcal{Q}$,
let $\hat{Q}^{s}_{t}$ denote the value of the QoI in scenario $s$ at
time $t$ given the predicted generation dispatches.
Eqs. \eqref{eq:risk_metrics:condition_expectation},
\eqref{eq:risk_metrics:reliability}, \eqref{eq:risk_metrics:risk} then
become
    \begin{align}
        \label{eq:risk_metrics:condition_expectation:proxy}
        \text{CE}(\mathcal{Q}_{t})
            &\approx
            \frac{
                \displaystyle \sum_{s \in \mathcal{S}} \hat{Q}^{s}_{t} \cdot \mathbbm{1}_{\hat{Q}^{s}_{t} \geq \hat{Q}_{t}(\alpha)}
            }{
                \displaystyle \sum_{s \in \mathcal{S}} \mathbbm{1}_{\hat{Q}^{s}_{t} \geq \hat{Q}_{t}(\alpha)}
            },\\
        \label{eq:risk_metrics:reliability:proxy}
        \mathbb{P}(\mathcal{Q}_{t} \geq \bar{Q}) 
            & \approx \frac{1}{S} \sum_{s \in \mathcal{S}} \mathbbm{1}_{\hat{Q}^{s}_{t} \geq \bar{Q}},\\
        \label{eq:risk_metrics:risk:proxy}
        \text{Risk}(\mathcal{Q}_{t}) 
            &= \mathbb{E}[c(\mathcal{Q}_{t})] \approx \frac{1}{S}\sum_{s \in \mathcal{S}} c(\hat{Q}^{s}_{t}),
    \end{align}
where $\hat{Q}_{t}(\alpha)$ denotes the $\alpha$ quantile of
$(\hat{Q}^{1}_{t}, ..., \hat{Q}^{S}_{t})$.

Observe that the proposed methodology predicts \emph{individual
  generation dispatches}, resulting in a risk assessment for {\em
  individual components} compared to existing coarser-grained
proposals, e.g., the zonal approach proposed in \cite{stover2022just}.

\begin{figure}[!t]
            \centering
            \includegraphics[width=0.9\columnwidth]{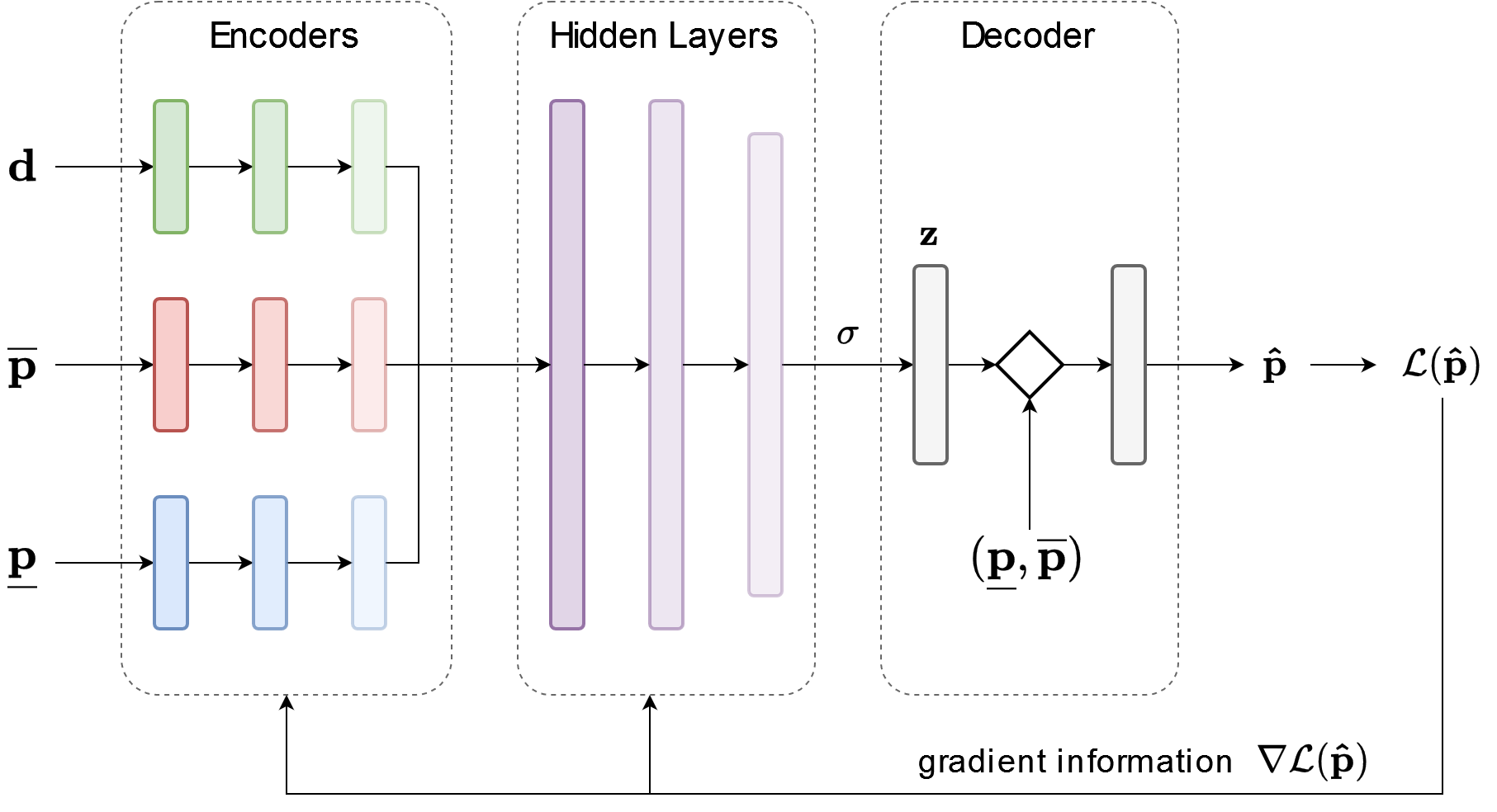}
            \caption{The Baseline DNN architecture. The hidden layers output $\z \, {\in} \, [0,1]^{N}$ using sigmoid activation $\sigma$, which is combined with min/max limits $\pgmin, \pgmax$ to ensure $\pgmin \, {\leq} \, \pghat \, {\leq} \, \pgmax$.
            The DeepOPF, DNN and E2ELR architectures extend the DNN architecture by adding a feasibility mechanism after the decoder.}
            \label{fig:ML:DNN}
        \vspace{-1em}
        \end{figure}

The proxies considered in this paper all extend the baseline DNN
architecture shown in Figure \ref{fig:ML:DNN}. They differ in how they
handle constraint violations. The baseline DNN consists of encoders,
hidden layers, and a decoder.  The DNN takes as input the load vector
$\pd$ and minimum/maximum limits $\pgmin, \pgmax$, which are each fed
to an encoder.  The encoders are small multi-layer perceptron (MLP)
models that use fully-connected layers with Rectified Linear Unit
(ReLU) activation
        \begin{align*}
            \text{ReLU}(x) = \max(0, x).
        \end{align*}
The output of the encoders, referred to as embeddings, are passed to
hidden layers, which are also fully connected layers with ReLU
activations.  A sigmoid activation is applied to the output of the
last hidden layer, which yields a vector $\z \, {\in} \, [0, 1]^{N}$.
The sigmoid function is defined as
        \begin{align}
            \sigma(x) = \frac{1}{1+e^{-x}} \in [0, 1].
        \end{align}
The predicted dispatch $\pghat$ is then obtained as 
        \begin{align}
            \pghat = \z \cdot \pgmin + (1-\z) \cdot \pgmax,
        \end{align}
where all operations are element-wise.  This approach ensures that
$\pgmin \, {\leq} \, \pghat \, {\leq} \, \pgmax$, which also ensures the feasibility
of the ramping constraints due to the structure of Algorithm
\ref{eq:risk_metrics:condition_expectation:proxy}. However, the
predicted dispatch may not satisfy power balance, i.e.,
$\e^{\top}\pghat = \e^{\top}\pd$ may not hold.

\subsection{The E2ELR Optimization Proxy}

The E2ELR architecture \cite{chen2023end} enforces min/max limits,
power balance constraints, and reserve constraints via dedicated,
differentiable feasibility layers. For the purpose of this paper, it
is sufficient to present the repair layer for power balance, which is
inserted between the approximation $\pghat$ and the computation of the
loss function. This repair layer takes $\pghat$ as input and outputs
        \begin{align}
            \label{eq:power_balance_layer}
            \begin{array}{rl}
               \pgfeas & = \left\{
               \begin{array}{ll}
                    (1 - \etaup) \pghat + \etaup \pgmax &  \text{if } \mathbf{e}^{\top} \pghat < \e^{\top} \pd\\
                    (1 - \etadn) \pghat + \etadn \pgmin &  \text{if } \mathbf{e}^{\top} \pghat \geq \e^{\top} \pd
               \end{array}
            \right.
            \end{array}
        \end{align}
        where $\pgmin \, {\leq} \, \pghat \, {\leq} \, \pgmax$, and $\etaup$, $\etadn$ are defined as follows:
        \begin{align}
            \label{eq:power_balance_layer:ratios}
            \etaup 
                &= \frac{
                    \mathbf{e}^{\top} \pd - \mathbf{e}^{\top} \pghat
                }{
                    \mathbf{e}^{\top} \pgmax - \mathbf{e}^{\top}\pghat
                },
            &
            \etadn 
                &= \frac{
                    \mathbf{e}^{\top} \pghat - \mathbf{e}^{\top}\pd
                }{
                    \mathbf{e}^{\top} \pghat - \mathbf{e}^{\top} \pgmin
                }.
        \end{align}
        The repaired dispatch $\pgfeas$ is feasible if and only if ED is feasible (see Theorem 1 in \cite{chen2023end}).

\subsection{The Training Methodology}
\label{sec:ML:training}

The training of E2ELR is based on self-supervised learning
\cite{park2023self,chen2023end}, which directly minimizes the
objective function \eqref{eq:ED:objective}.  Consider a dataset with
$M$ data points
\begin{align}
        \mathcal{D} =\Bigl\{ \x^1, \x^2, \cdots, \x^M \Bigl \},
\end{align}
where $\x^t = \{\pd^t, \pgmax^t, \pgmin^t \}$, Let $f_{\btheta}$ be
the machine learning model with trainable parameters $\btheta$.
The training of $f_{\btheta}$ amounts to solving the optimization
problem
\begin{align}
\btheta^* = \argmin_{\btheta} \frac{1}{M} \sum^{M}_{t=1} \mathcal{L}(f_{\btheta}(\x^t)),
\end{align}
In E2ELR, the loss function is the objective function of ED:
\begin{align}
\mathcal{L}(\pghat) &= c(\pghat) + \Mth \|\xith(\pghat)\|_1.
\end{align}
Observe that E2ELR only uses a set of inputs for training: it does not
require labeled data.

\section{Experimental Results}
\label{sec:results}

\newcommand{\ieeeSmall}{\texttt{ieee300}}
\newcommand{\pegase}{\texttt{pegase1k}}

\subsection{The Machine Learning Architectures}
\label{sec:methdology:ML}

A key advantage of \raproxy{} is its ability to produce ED-feasible
solutions quickly at each stage of the MC simulations. To demonstrate
these benefits, the paper compares \raproxy{} to other ML
architectures, which offer no feasibility guarantees. To train these
infrastructures, the loss function should be generalized to
$\mathcal{L}(\pghat) = \varphi(\pghat) + \lambda \psi(\pghat)$, where the first term
$\varphi(\pghat)$ is the objective value of ED
    \begin{align}
        \label{eq:SSL:objective}
        \varphi(\pghat) &= c(\pghat) + \Mth \|\xith(\pghat)\|_1
    \end{align}
    and the second term $\psi(\pghat)$ penalizes constraint violations. For instance,
    the power balance violations can be formalized as
 \begin{align}
 \label{eq:train:constraint_penalty}
 \psi(\pghat) = \Mpb |\mathbf{e}^{\top}\pd - \mathbf{e}^{\top} \pghat |,
 \end{align}
 where $\Mpb$ are penalty coefficients.

The evaluation considers three ML architectures to show the benefits
of \raproxy{}: the baseline DNN architecture presented earlier,
Deep-OPF \cite{pan2020deepopf}, and DC3 \cite{donti2021dc3}.

\subsubsection{DeepOPF}
\label{sec:ML:DeepOPF}
    
The DeepOPF architecture \cite{pan2020deepopf} extends the DNN model
by adding a equality completion step that enforces system power balance
\eqref{eq:ED:power_balance}. Let $\pghat$ denote the output of DNN and
assume bus $1$ is the slack bus without loss of generality. The
equality completion layer outputs $\pgfeas \, {\in} \,
\mathbb{R}^{N}$, where
\begin{align}
\label{eq:equality_elimination}
\pgfeas_{1} &= \e^{\top} \pd - \sum_{i > 1} \pghat_{i}, \ 
&
\pgfeas_{i} &= \pghat_{i}, \ \forall i > 1,
\end{align}
thus ensuring the satisfaction of the power balance constraint
\eqref{eq:ED:power_balance}. However, the dispatch of the slack
generator $\pgfeas_{1}$ is not guaranteed to respect its minimum and
maximum limits.

\subsubsection{Deep Constraint Completion and Correction (DC3)}
\label{sec:ML:DC3}

 The DC3 architecture \cite{donti2021dc3} uses the same equality
 completion mechanism as DeepOPF, with an additional inequality
 correction step.  The latter mitigates potential violations of the
 min/max limits for the slack generator, and uses unrolled gradient
 descent mechanism that minimizes the violation of inequality
 constraints.  The number of gradient steps in the correction
 mechanism is set to a fixed value $K$.  If $K$ is large enough, then
 DC3 returns feasible solutions, albeit at the cost of longer
 inference times.  A smaller $K$ improves inference time, but may
 yield constraint violations.

\subsection{Data Generation}
\label{sec:results:data}

Experiments are carried on the IEEE 300-bus system (\ieeeSmall)
\cite{PSTCA} and the Pegase 1354-bus system (\pegase)
\cite{josz2016_Pegase}, using data from the UnitCommitment (UC)
package \cite{Xavier2022_UnitCommitment_jl}.  The experiments use the
network topology and cost information in the original test cases.
Load time series and the minimum/maximum limits and ramping rates of
the generators are obtained from \cite{Xavier2022_UnitCommitment_jl}.
For ease of presentation, all generators are assumed to be online; the
proposed framework naturally extends to varying commitment decisions.
For each system, \cite{Xavier2022_UnitCommitment_jl} provides 365
benchmark instances, one per day in 2017, with each instance spanning
a horizon of $T \, {=} \, 36$ hours.  This initial dataset of 365 load
profiles is augmented by generating addition load samples, using the
methodology proposed in \cite{Xavier2021_LearningToSolveSCUC}.

For a given system and a given day, bus-level load scenarios are of the form
$
    \pd_{i,t} = \xi_{i,t} \times (\gamma_{t} \times \bar{P}^{d}),
$
where $\bar{P}^{d}$ denotes the maximum total load in the reference UC
instance, $\gamma \, {\in} \, \mathbb{R}^{T}$ captures the
distribution of load-to-peak ratio, and $\xi_{i,t}$ denotes bus-level
scaling factors.  The distribution of $\gamma$ is estimated from the
original 365 instances using a Gaussian copula model with Beta
marginal distributions.  Note that, in the instances provided by
\cite{Xavier2022_UnitCommitment_jl}, each bus $i$ always accounts for
a fixed proportion $\xi^{\text{ref}}_{i}$ of total demand.  Therefore,
nodal disaggregation factors $\xi$ are of the form $\xi_{i,t} \, {=}
\, \xi^{\text{ref}}_{i} \, {\times} \, \eta_{i, t}$, where
$\eta_{i,t}$ is uncorrelated white noise with log-normal distribution
of mean $1$ and 5\% standard deviation.

\begin{figure}[!t]
    \centering
    \includegraphics[width=0.44\columnwidth]{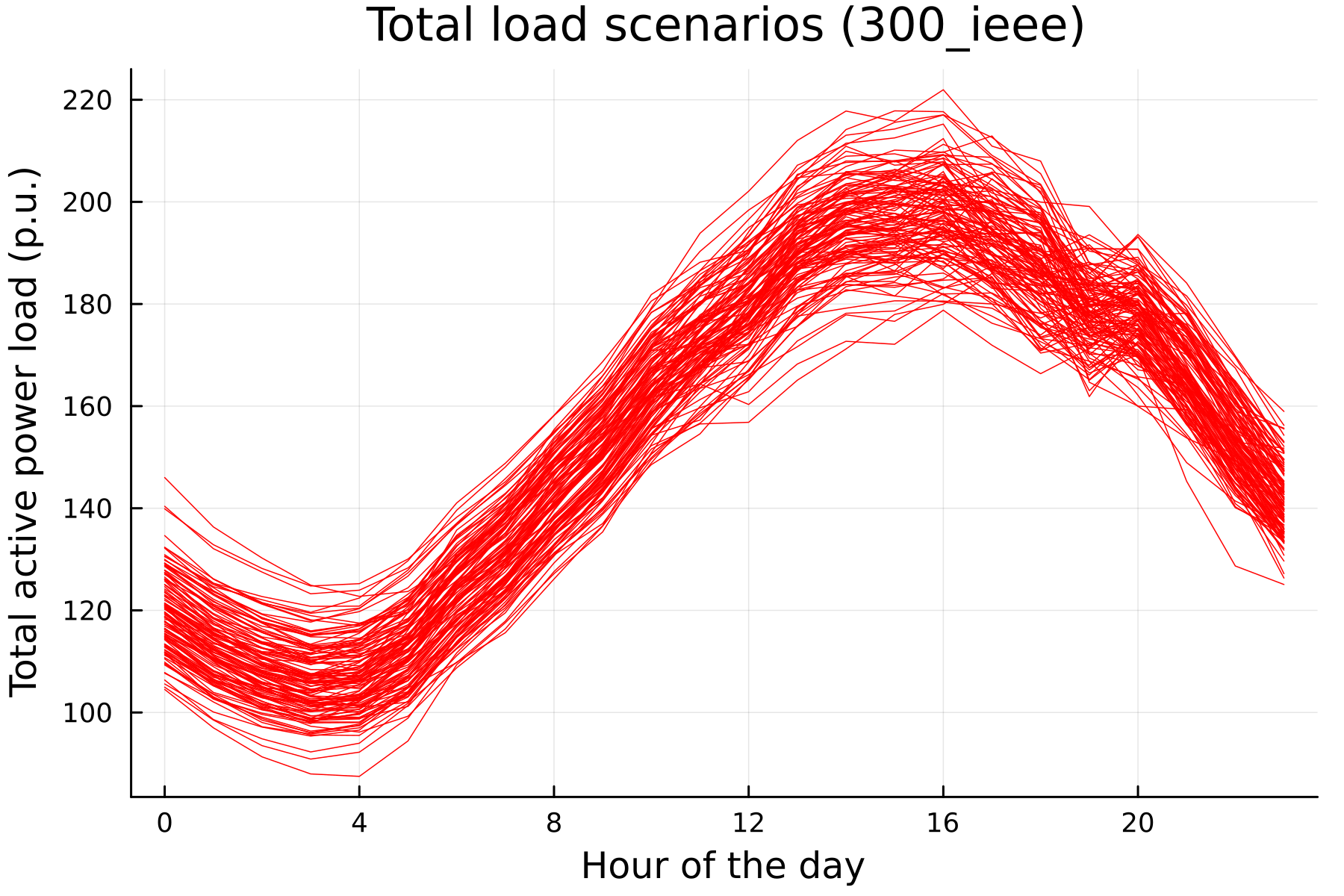}
    \hfill
    \includegraphics[width=0.44\columnwidth]{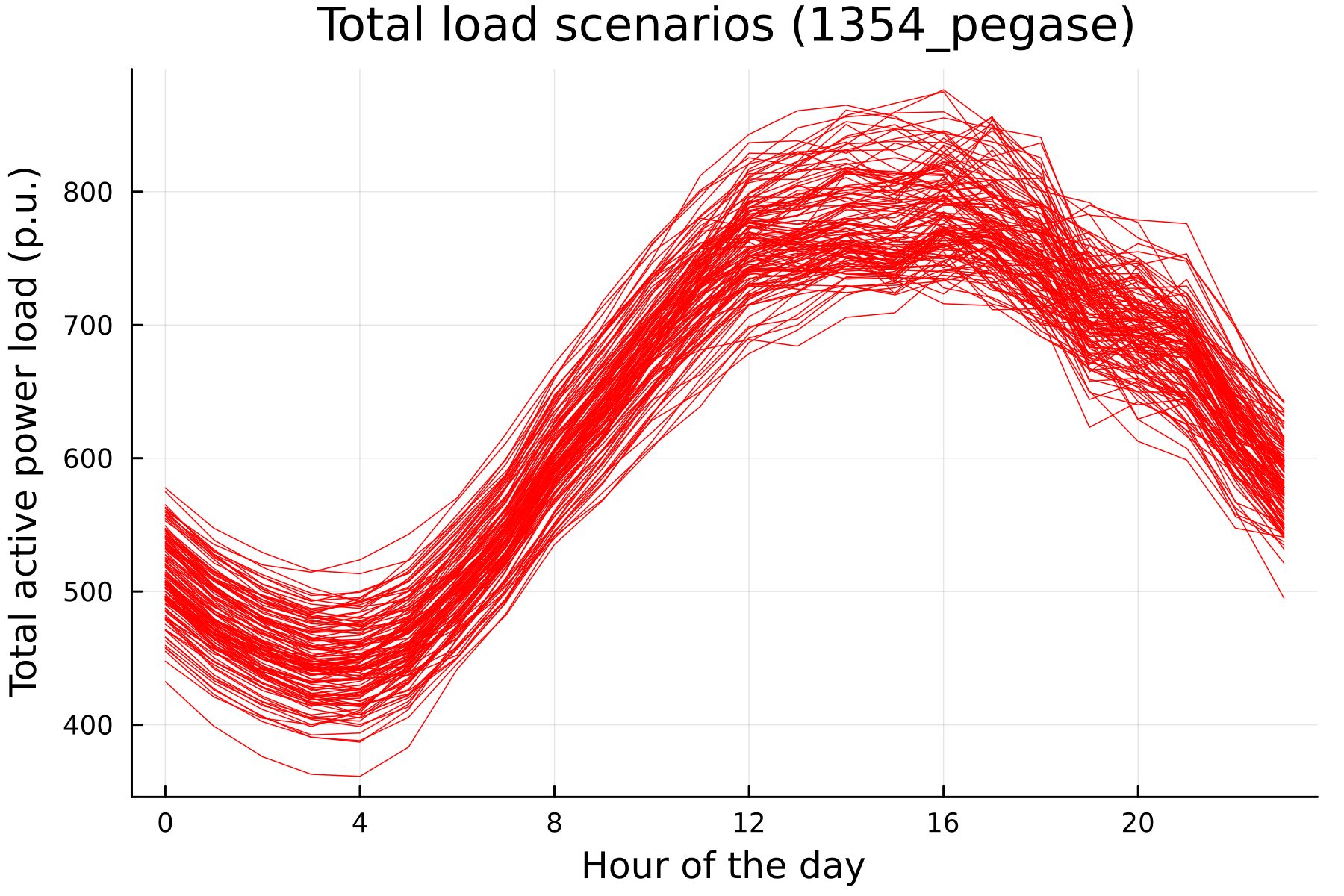}
    \caption{Total load scenarios used in the experiments. Left: \ieeeSmall{}, on August 22nd, 2017. Right: \pegase{}, on July 13th, 2017.}
    \label{fig:load_scenarios}
    \vspace{-1em}
\end{figure}

Figure \ref{fig:load_scenarios} depicts daily total load scenarios for
the \ieeeSmall{} and \pegase{} systems.  The paper selected August
22nd, 2017 for \ieeeSmall, which was found to encounter capacity
shortages, and July 13th, 2017 for \pegase{}, which is the day that
experienced the highest load over the year.  A total of 1024 and 2048
scenarios are generated for the 300 and 1354-bus system, respectively.
Each dataset is then split between training (80\%), validation (10\%)
and testing (10\%).

\subsection{Implementation Details}
\label{sec:results:implementation}

    The optimization problems are formulated in Julia using JuMP \cite{Dunning2017_JuMP}, and solved with Gurobi 9.5 \cite{gurobi}, with a single CPU thread and its default parameters.
    All deep learning models are implemented using PyTorch \cite{paszke2017automatic} and trained using the Adam optimizer \cite{kingma2014adam}.
    All models are hyperparameter tuned using a grid search, with the learning rate from $\{10^{-2}, 10^{-3}\}$, hidden dimension from $\{128, 256\}$ and penalty term of constraint violation from $\{0.1, 1\}$.  
    For each system, the model with the lowest validation loss is selected and the performances on the test set are reported.
    During training, the learning rate is reduced by a factor 10 if the validation loss shows no improvement for a consecutive sequence of 10 epochs. 
    In addition, training is stopped if the validation loss does not improve for consecutive 20 epochs. 
    Experiments are conducted on dual Intel Xeon 6226@2.7GHz machines running Linux, on the PACE Phoenix cluster \cite{PACE}.  
    The ML models are trained on Tesla V100-PCIE GPUs with 16GBs HBM2 RAM.

\subsection{Comparison Framework}
\label{sec:results:comparison}

The risk simulations follow Algorithm \ref{alg:risk}.  Since all
simulations start from the same initial generation setpoint, ramping
issues may arise at the beginning of the simulation, which results in
artificially high risk during the first few hours of the day.
Therefore, results are reported for a 24-hour period starting at 4am
and ending at 4am.  Finally, to ensure fair comparison between all ML
models, the paper implements the following mechanism.  All generation
dispatches are clipped, if necessary, to ensure that minimum/maximum
limits \eqref{eq:ED:bounds} are satisfied.  This mechanism is executed
at each time step of every simulation, and captures the fact that
generators cannot adjust their output beyond their physical
capabilities.  While this does not affect the output of DNN and E2ELR,
it may modify the output of DeepOPF and DC3, thus resulting in power
balance violations instead of violations of the slack bus constraints.

% \subsection{Performance Metrics}
% \label{sec:results:metrics}

%     This work considers the operation cost, power balance violation, total thermal violation, and thermal violation on individual lines as the quality of interests.
%     For the operation cost, this paper reports scatter plots and Q-Q plots to show the distribution of the costs generated using proxies and the optimization solver.
%     For the others, the conditional expectation, reliability, and risk are reported.
    % Highlight that the thermal violation on individual lines could be simulated because of the high fidelity of the optimization proxies. 

% This work considers the operation cost, power balance violation, and thermal violation as the quality of interests.
% For the operation cost, this paper reports Q-Q plots, and CVaR to show the distribution of the costs generated using proxies and the optimization solver.
% For the others, the CVaR, probability of failure, and risk are reported.

\subsection{Global Risk Analysis}

\begin{figure}[!t]
        \centering
        \subfloat[
            Level-1 (CVaR) risk metrics for \ieeeSmall{}.
            Left: $\text{CVaR}_{90\%}$ of system power imbalance.
            Right: $\text{CVaR}_{90\%}$ of total thermal violations.
        ]{
            \label{fig:risk:300_ieee:L1}
            \includegraphics[width=0.44\columnwidth]{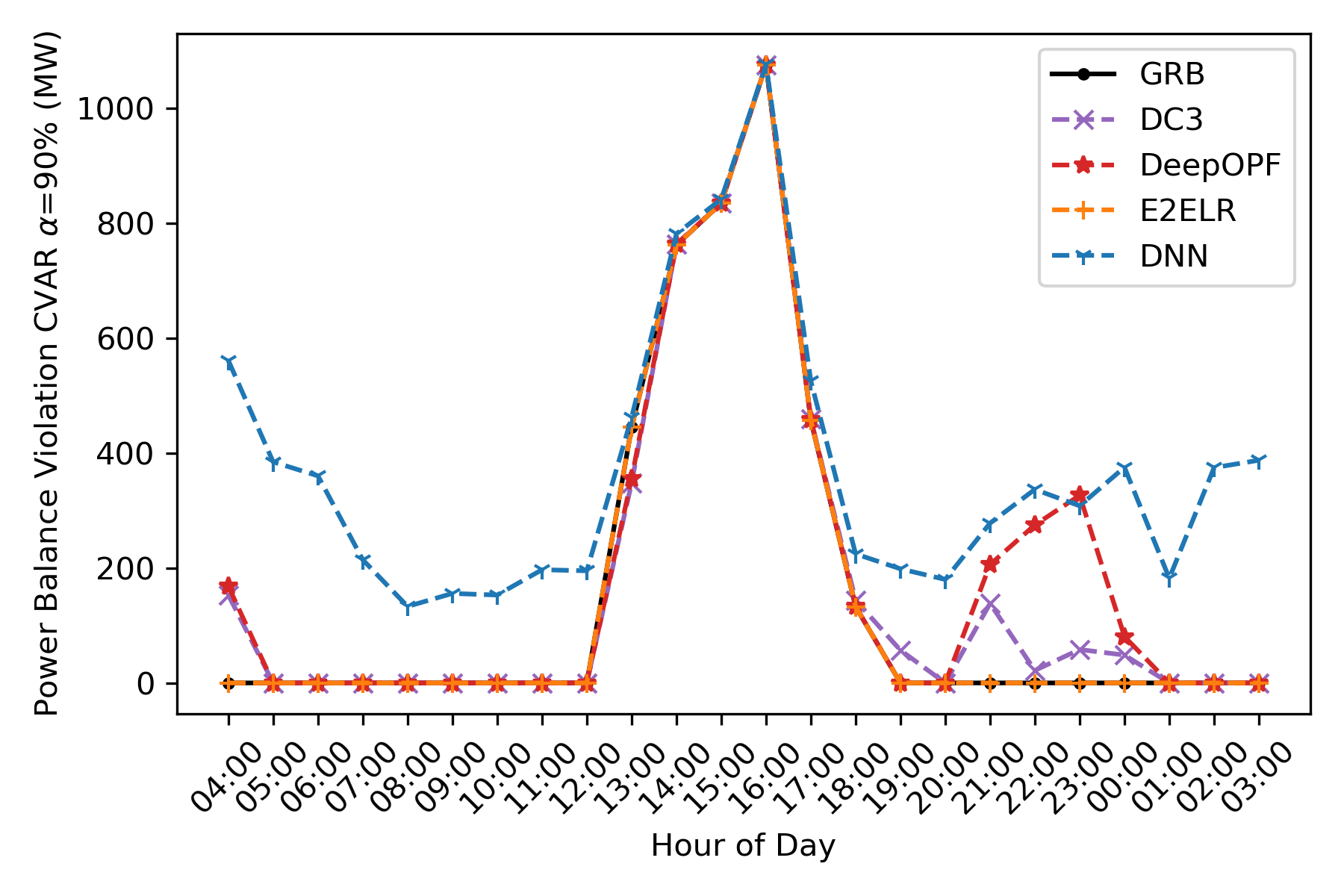}
            \hspace{1em}
            \includegraphics[width=0.44\columnwidth]{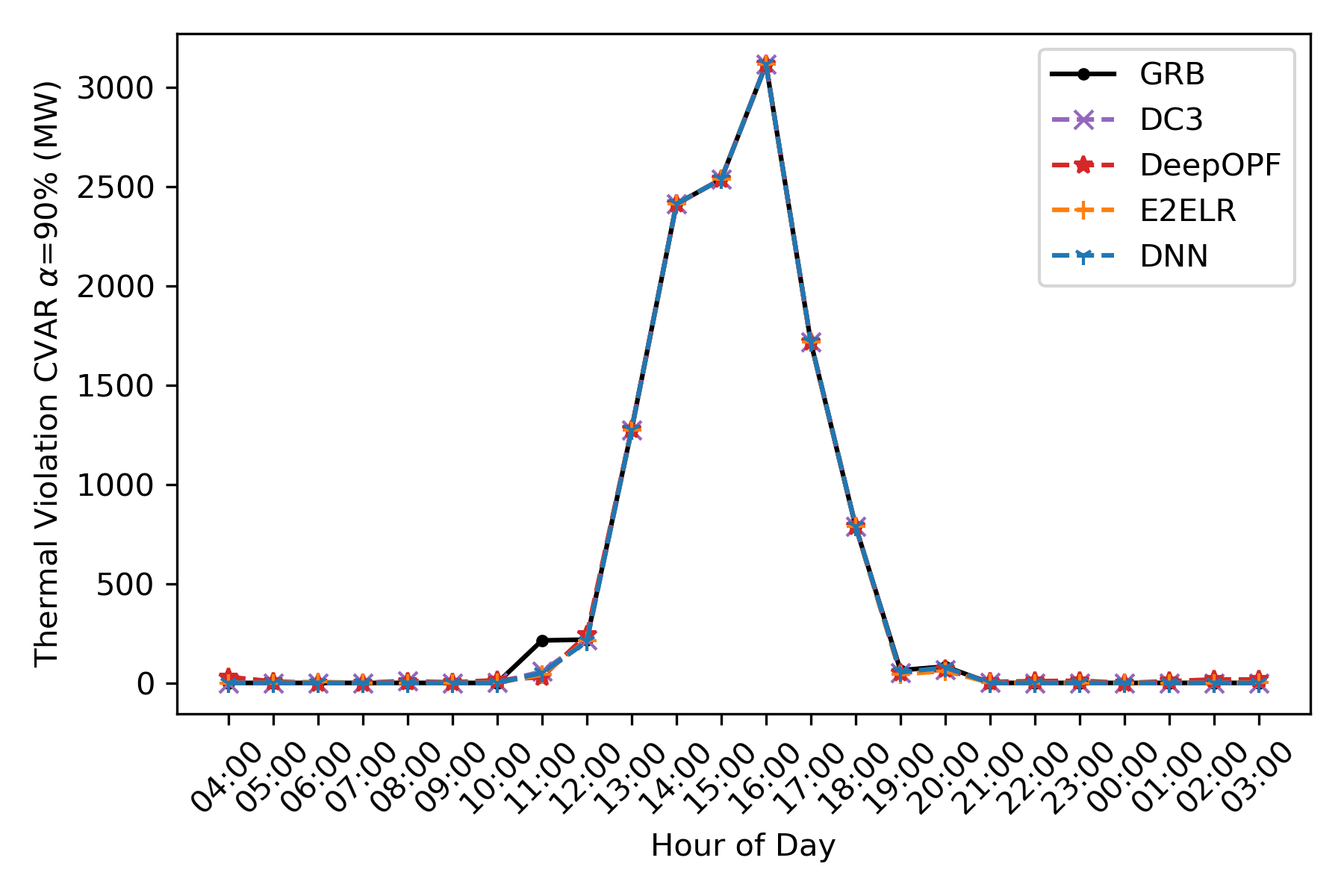}
        }\\
        \subfloat[
            Level-2 (probability of adverse event) risk metrics for \ieeeSmall{}.
            Left: probability of system power imbalance.
            Right: probability of non-zero total thermal violations.
        ]{
            \label{fig:risk:300_ieee:L2}
            \includegraphics[width=0.44\columnwidth]{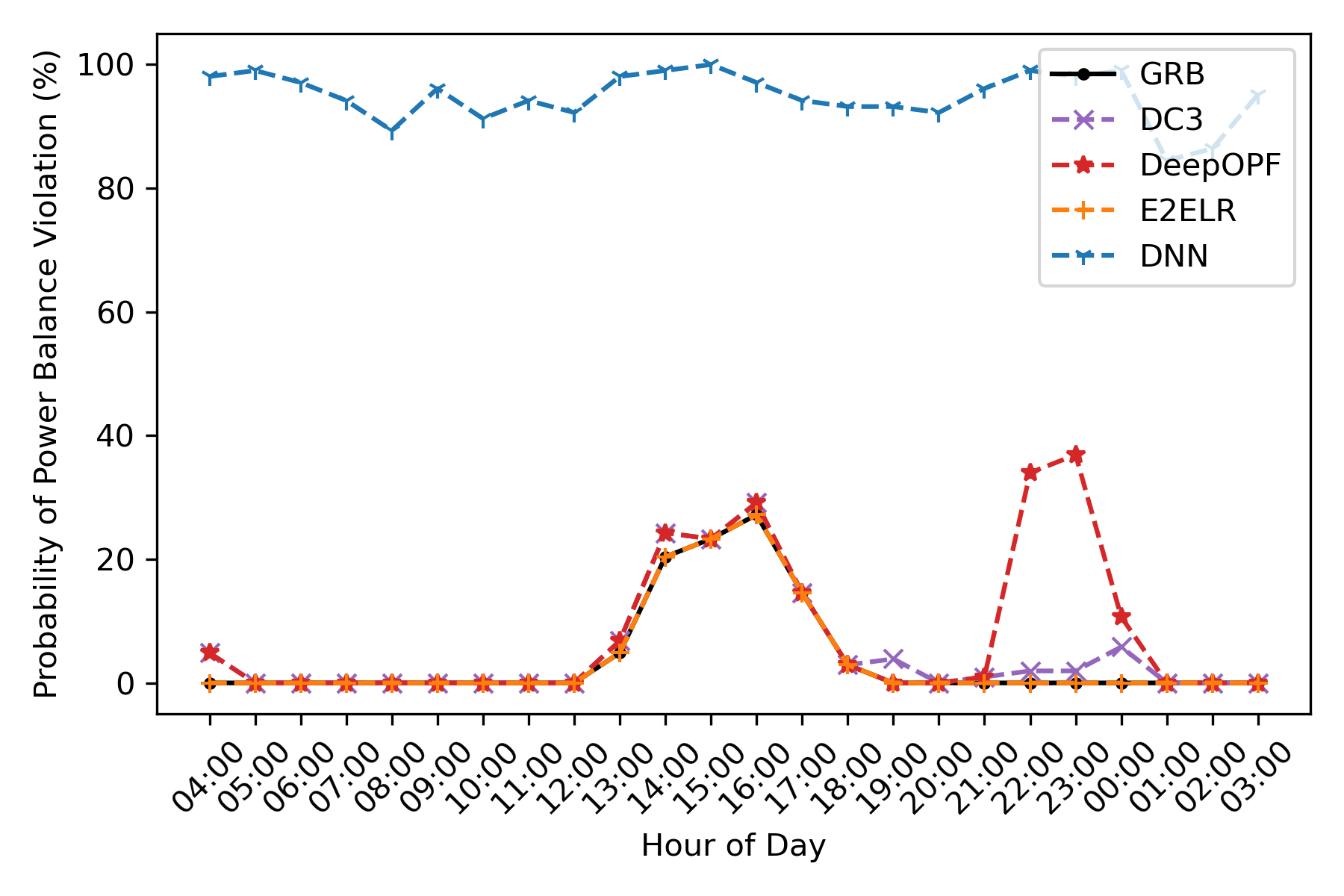}
            \hspace{1em}
            \includegraphics[width=0.44\columnwidth]{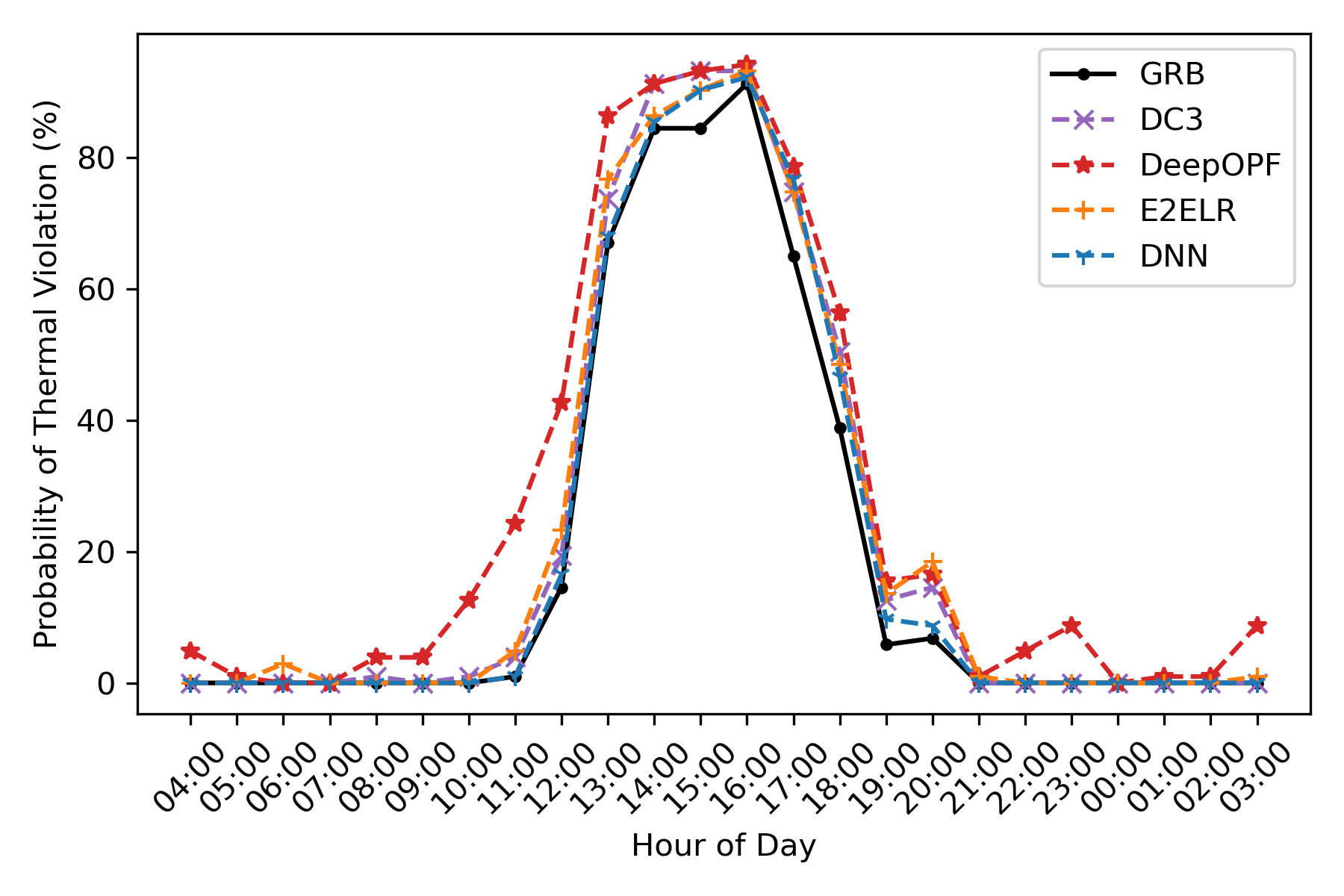}
        }\\
        \subfloat[
            Level-3 (financial risk) metrics for \ieeeSmall{}.
            Left: risk of system power imbalance, in log scale.
            Right: risk of total thermal violations.
        ]{
            \label{fig:risk:300_ieee:L3}
            \includegraphics[width=0.44\columnwidth]{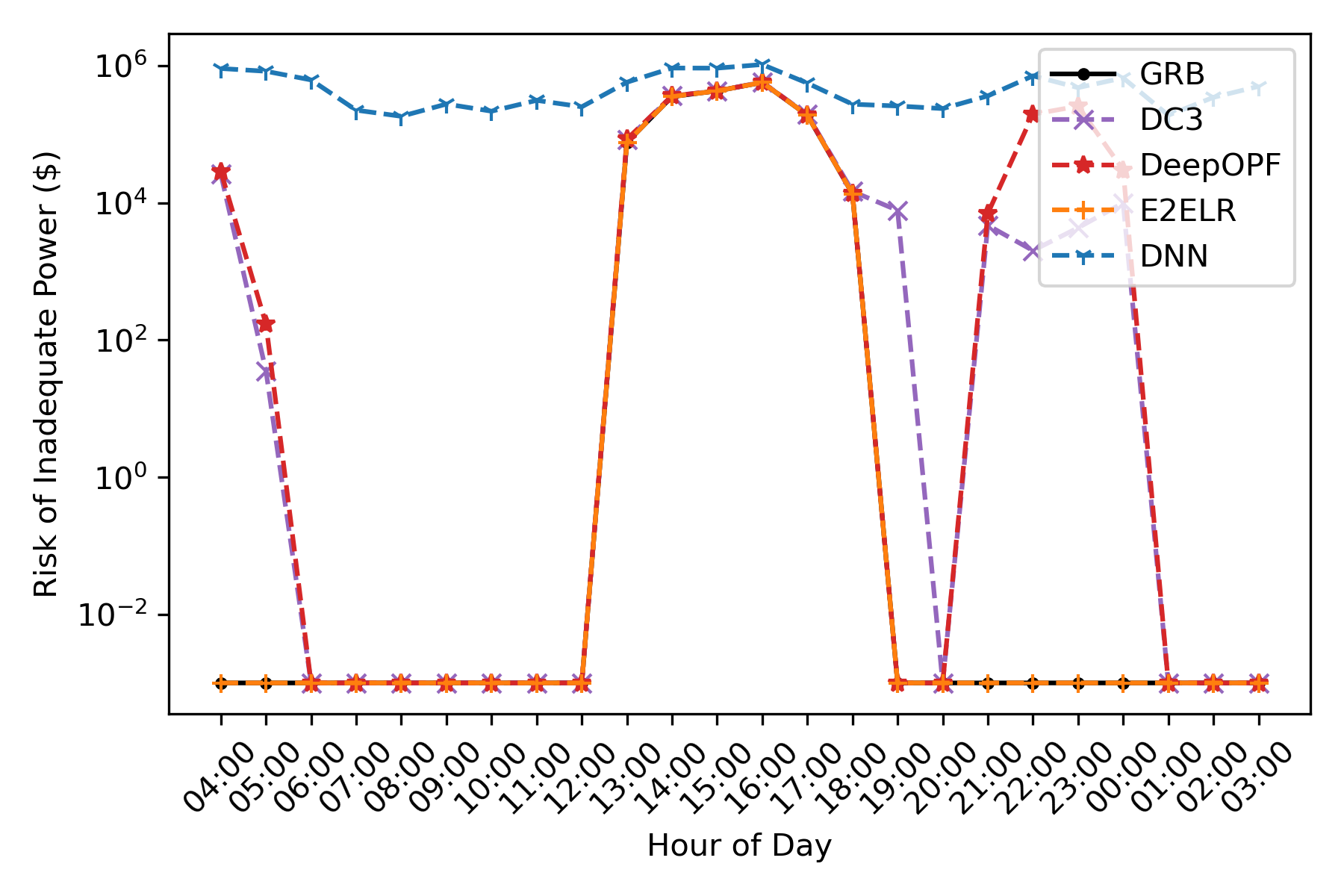}
            \hspace{1em}
            \includegraphics[width=0.44\columnwidth]{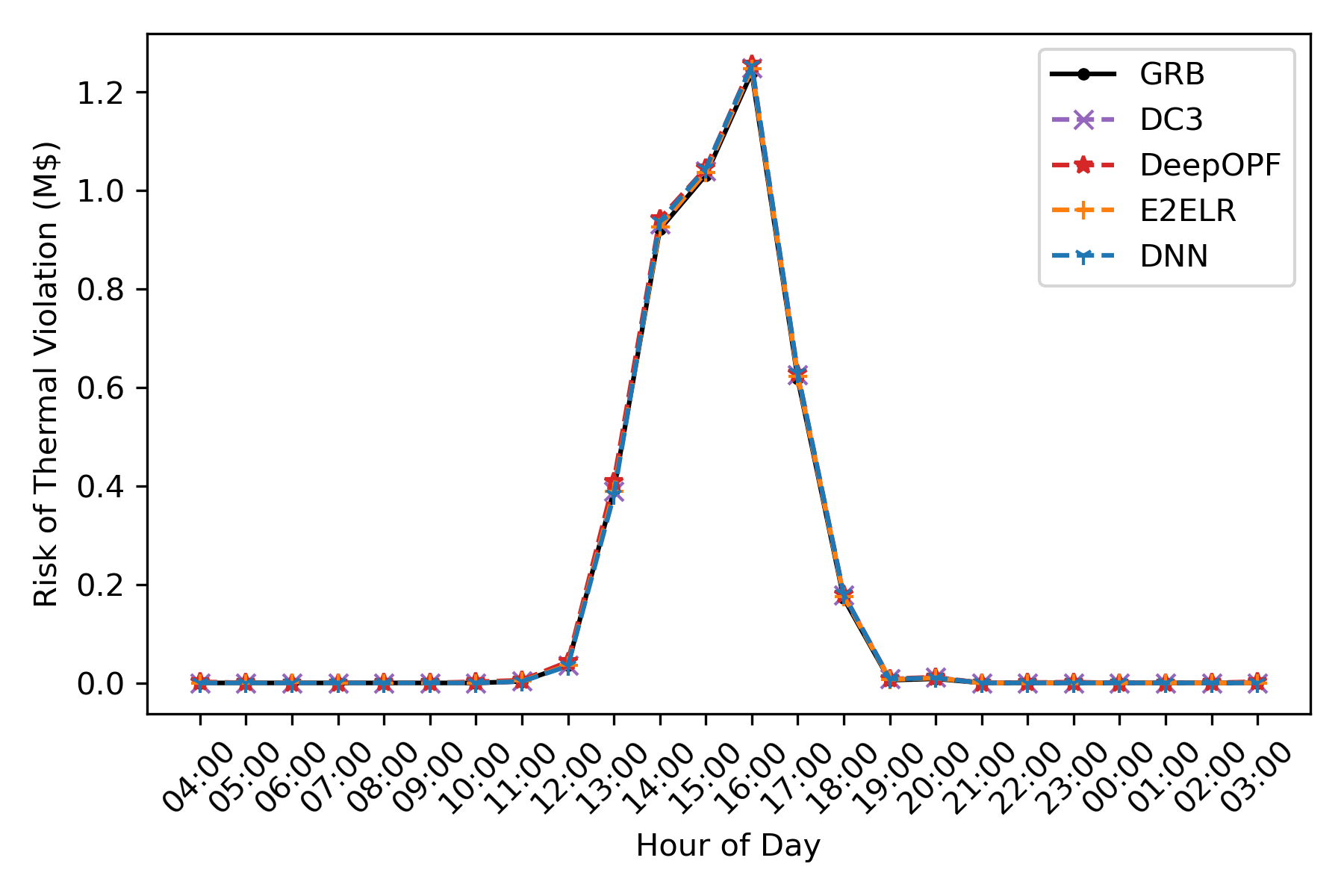}
        }\\
        \caption{Risk assessment results on the \ieeeSmall{} system.}
        \label{fig:risk:300_ieee}
        \vspace{-1em}
    \end{figure}

\begin{figure}[!t]
        \centering
        \subfloat[
            Level-1 (CVaR) risk metrics for \pegase{}.
            Left: $\text{CVaR}_{90\%}$ of system power imbalance.
            Right: $\text{CVaR}_{90\%}$ of total thermal violations.
        ]{
            \label{fig:risk:1354_pegase:L1}
            \includegraphics[width=0.44\columnwidth]{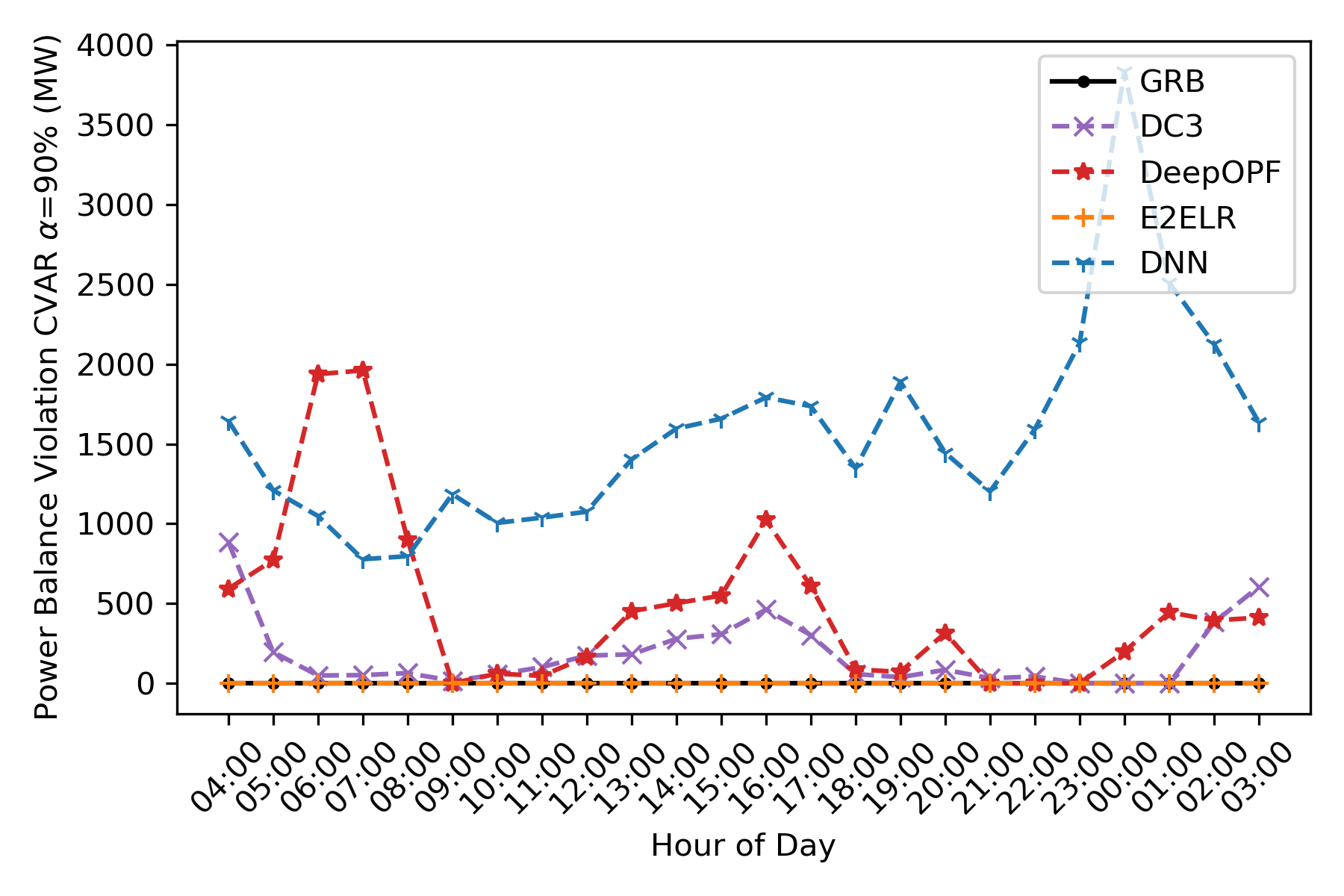}
            \hspace{1em}
            \includegraphics[width=0.44\columnwidth]{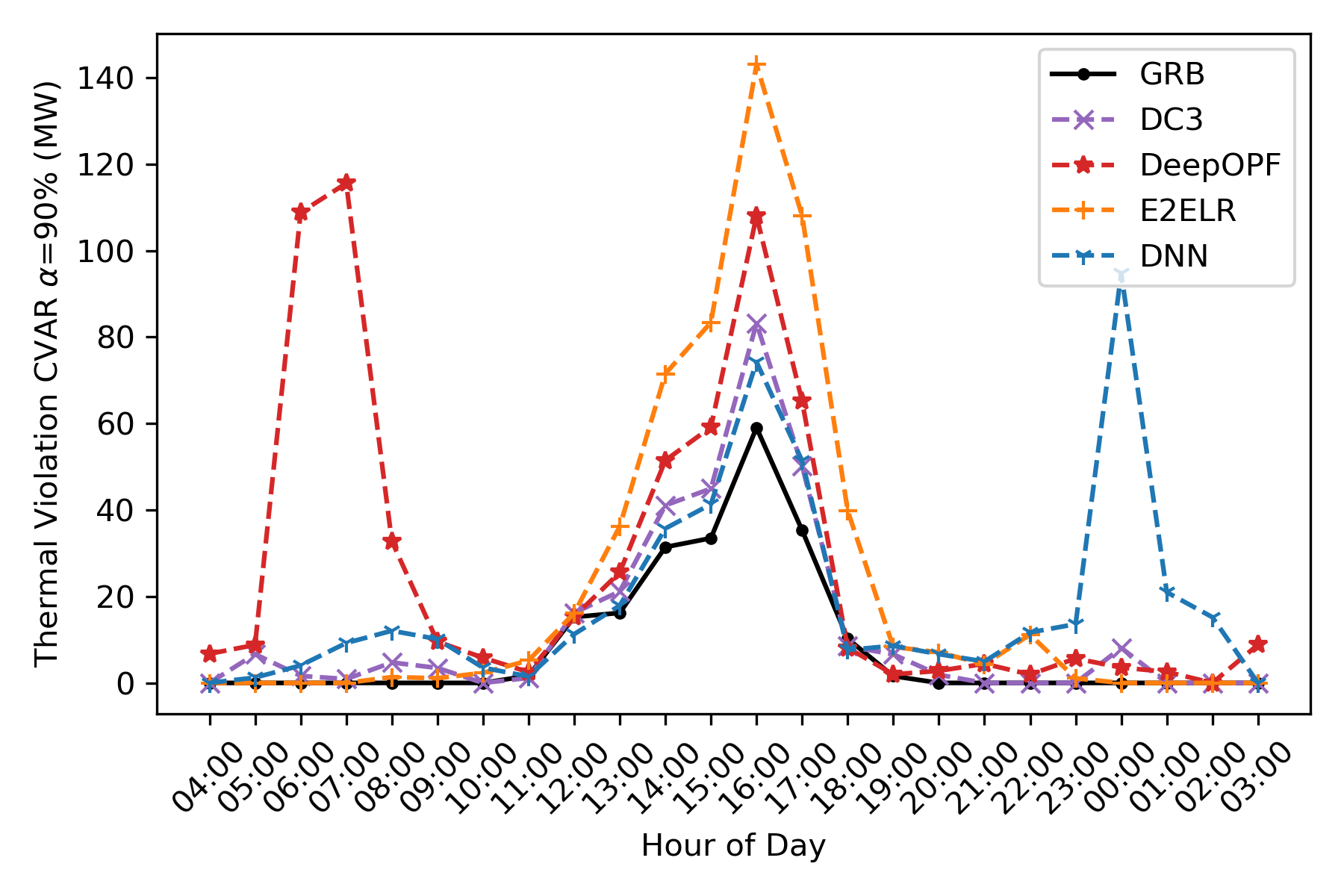}
        }\\
        \subfloat[
            Level-2 (probability of adverse event) risk metrics for \pegase{}.
            Left: probability of system power imbalance.
            Right: probability of non-zero total thermal violations.
        ]{
            \label{fig:risk:1354_pegase:L2}
            \includegraphics[width=0.44\columnwidth]{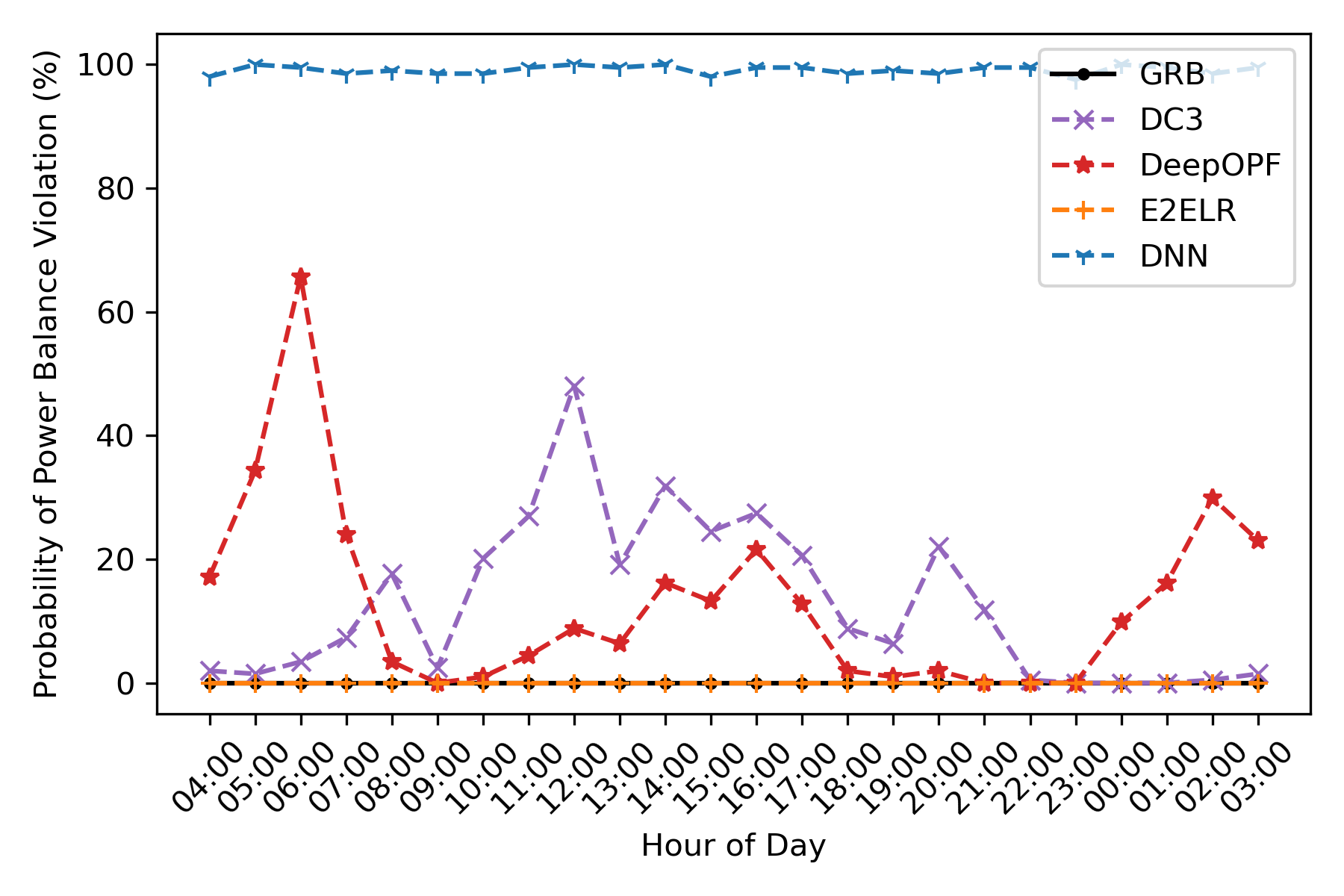}
            \hspace{1em}
            \includegraphics[width=0.44\columnwidth]{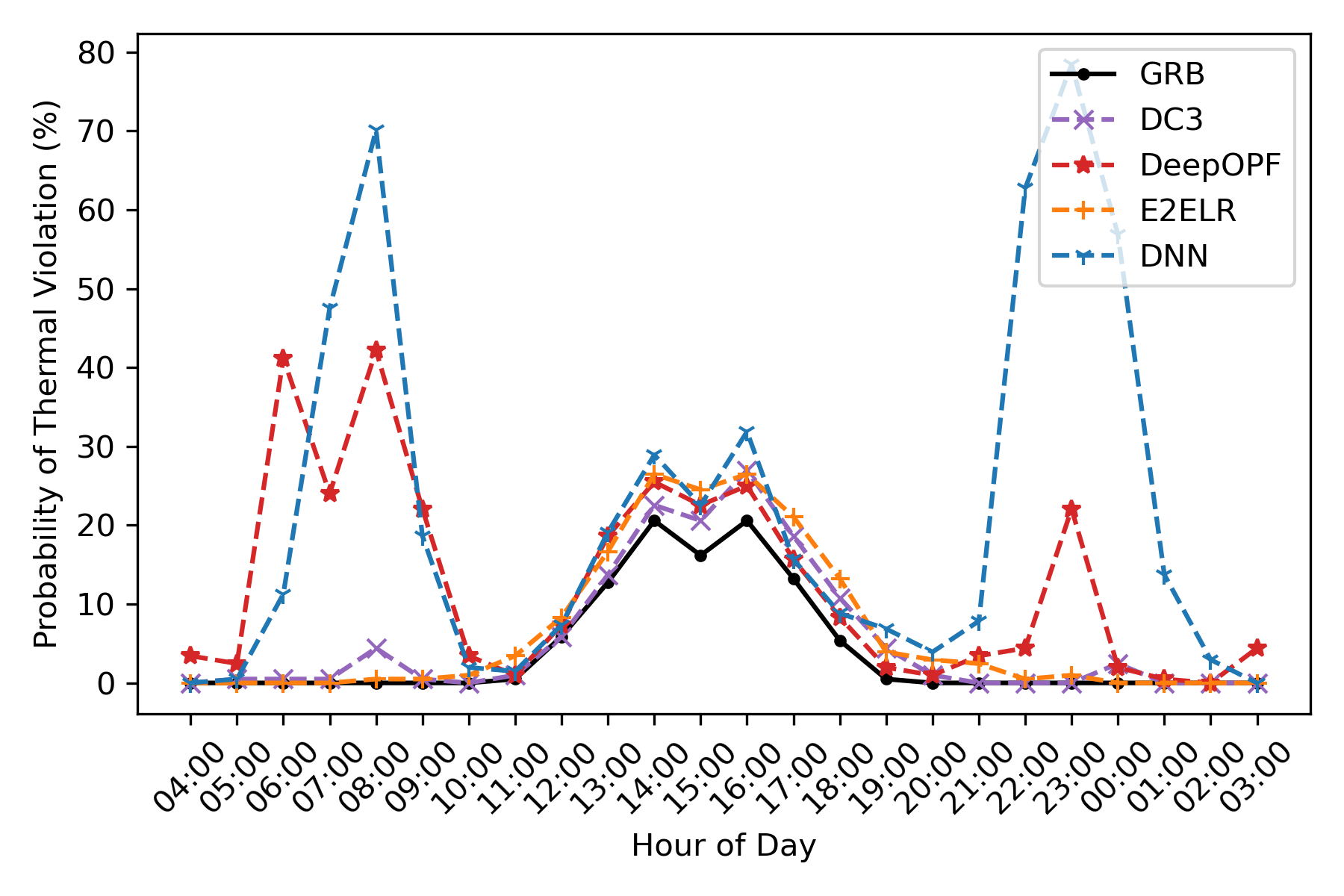}
        }\\
        \subfloat[
            Level-3 (financial risk) metrics for \pegase{}.
            Left: risk of system power imbalance, in log scale.
            Right: risk of total thermal violations.
        ]{
            \label{fig:risk:1354_pegase:L3}
            \includegraphics[width=0.44\columnwidth]{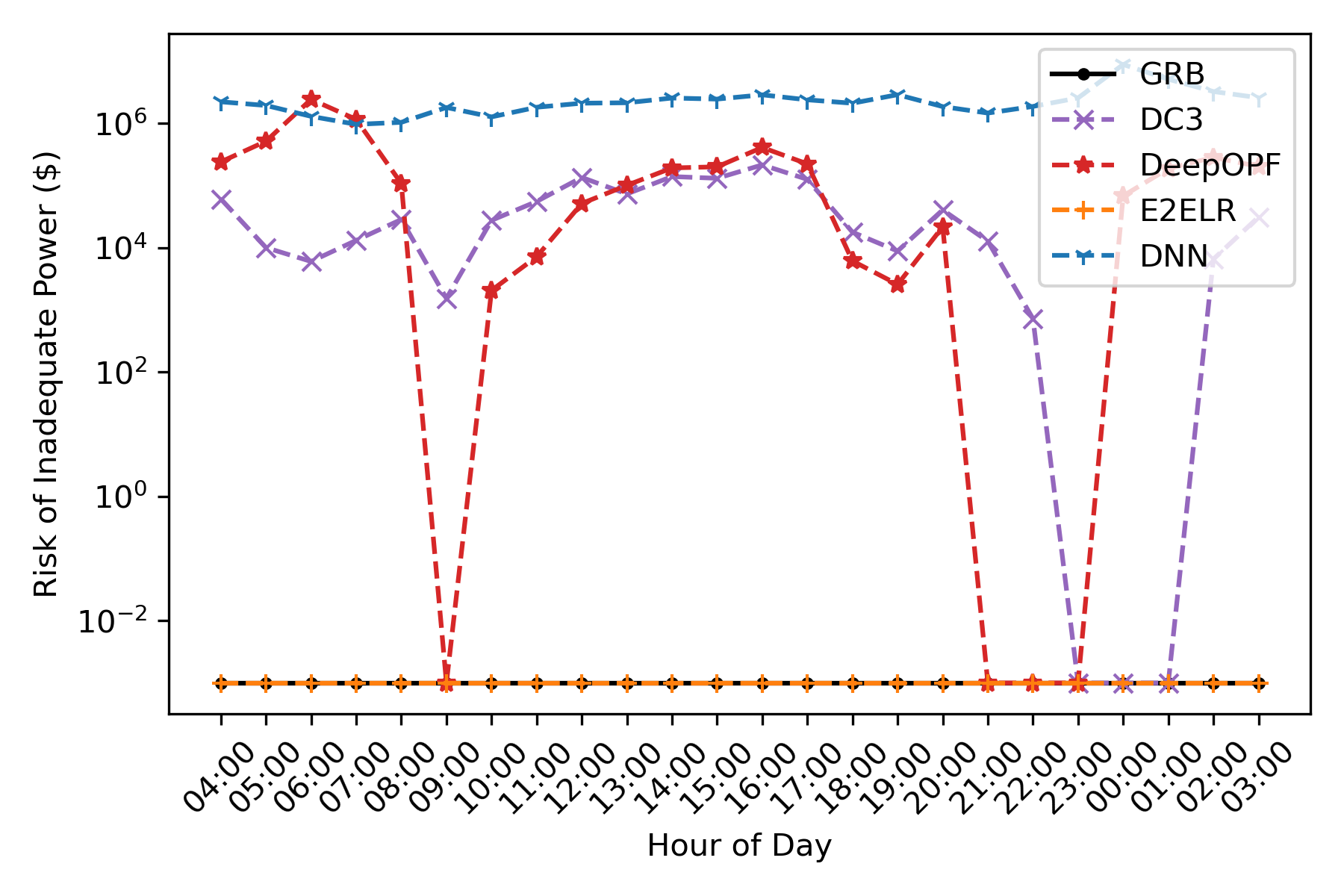}
            \hspace{1em}
            \includegraphics[width=0.44\columnwidth]{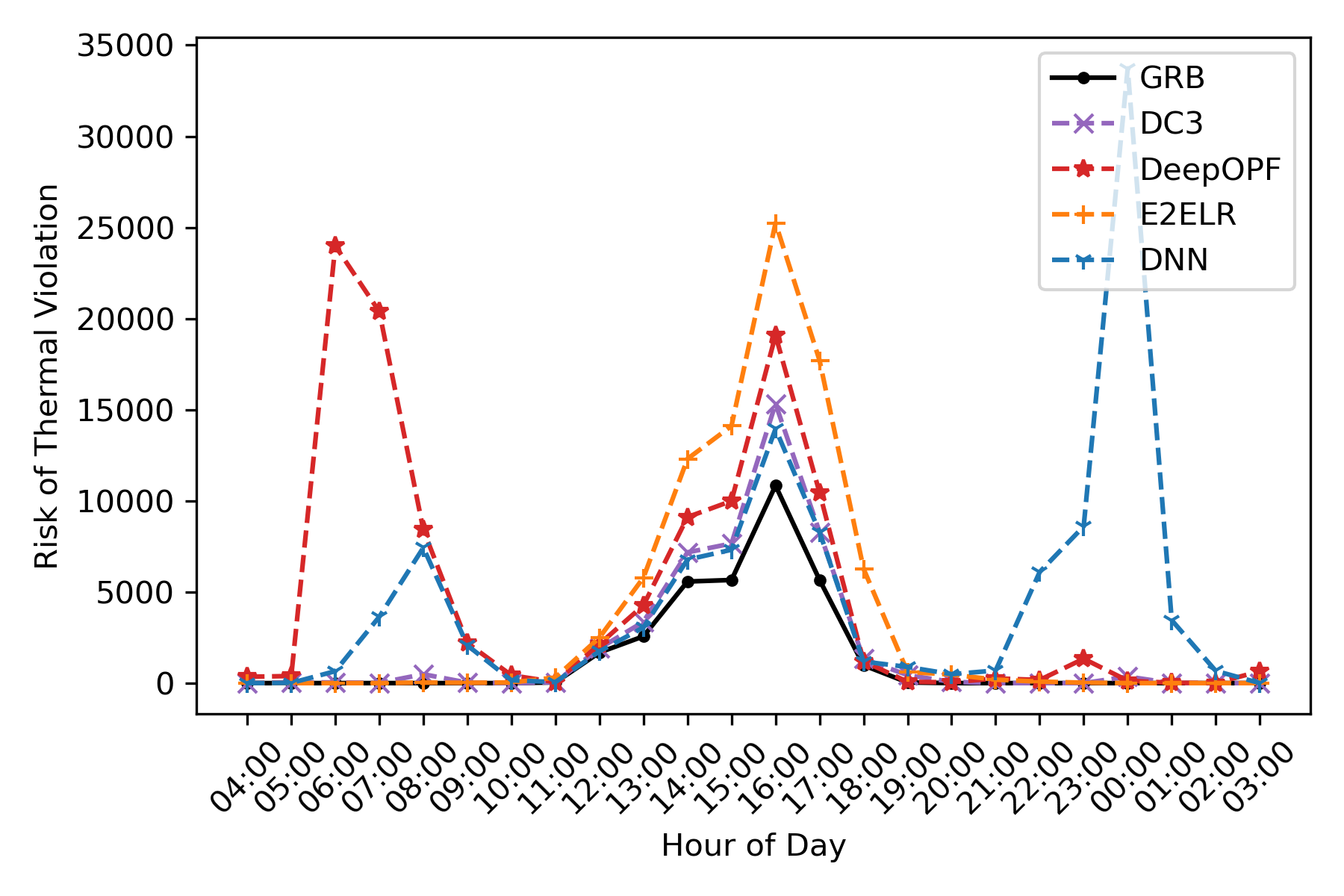}
        }\\
        \caption{Risk assessment results on the \pegase{} system.}
        \vspace{-1em}
        \label{fig:risk:1354_pegase}
    \end{figure}

Figures \ref{fig:risk:300_ieee} and \ref{fig:risk:1354_pegase} compare
ground truth and proxy-based risk profiles for the \ieeeSmall{} and
\pegase{} systems, respectively.  The ground truth risk profiles are
obtained by executing Algorithm \ref{alg:risk} with an optimization
solver (GRB).  The proxy-based risk profiles are obtained by replacing
the optimization solver with an optimization proxy.  Each figure
reports, for each hour: the 90\%-CVaR (level 1) for system-wide power
imbalance and total thermal violations, the probability of system-wide
power imbalance and thermal violations, and the risk of power
imbalance and total thermal violations.

First observe that RA-DNN, RA-DC3 and RA-DeepOPF consistently over-estimates
power imbalances.  This is most striking for RA-DNN, which almost always
predicts non-zero power imbalances, while RA-DeepOPF and RA-DC3 wrongly
predict system imbalances during the afternoon hours on \ieeeSmall{}
and throughout the day on \pegase.  This is because the DNN, DeepOPF,
and DC3 are not guaranteed to output solutions that satisfy all ED
constraints.  In contrast, the E2ELR architecture, which always
outputs ED-feasible solutions, prefectly matches the ground truth risk
assessment for power imbalances.
    
Second, all proxies tend to over-estimate thermal violations, which is
more flagrant on the larger \pegase{} system (Figure
\ref{fig:risk:1354_pegase}). RA-DNN, RA-DeepOPF and, to a lesser extent,
RA-DC3, incorrectly predict thermal violations during the morning hours
(5am--9am) and during the night (10pm--2am).  While the E2ELR-based
risk profile over-estimates thermal violations during the afternoon
(12pm--5pm), Figure \ref{fig:risk:1354_pegase:L2} shows that RA-DC3 and
RA-E2ELR correctly identify those hours during which congestion occurs.
It is important to note that inaccurate risk predictions may result in
unnecessary preventive actions, which result in higher economic cost
for the system.

\subsection{Global Cost Analysis}

\begin{figure}[!t]
        \centering
        \includegraphics[width=0.44\columnwidth]{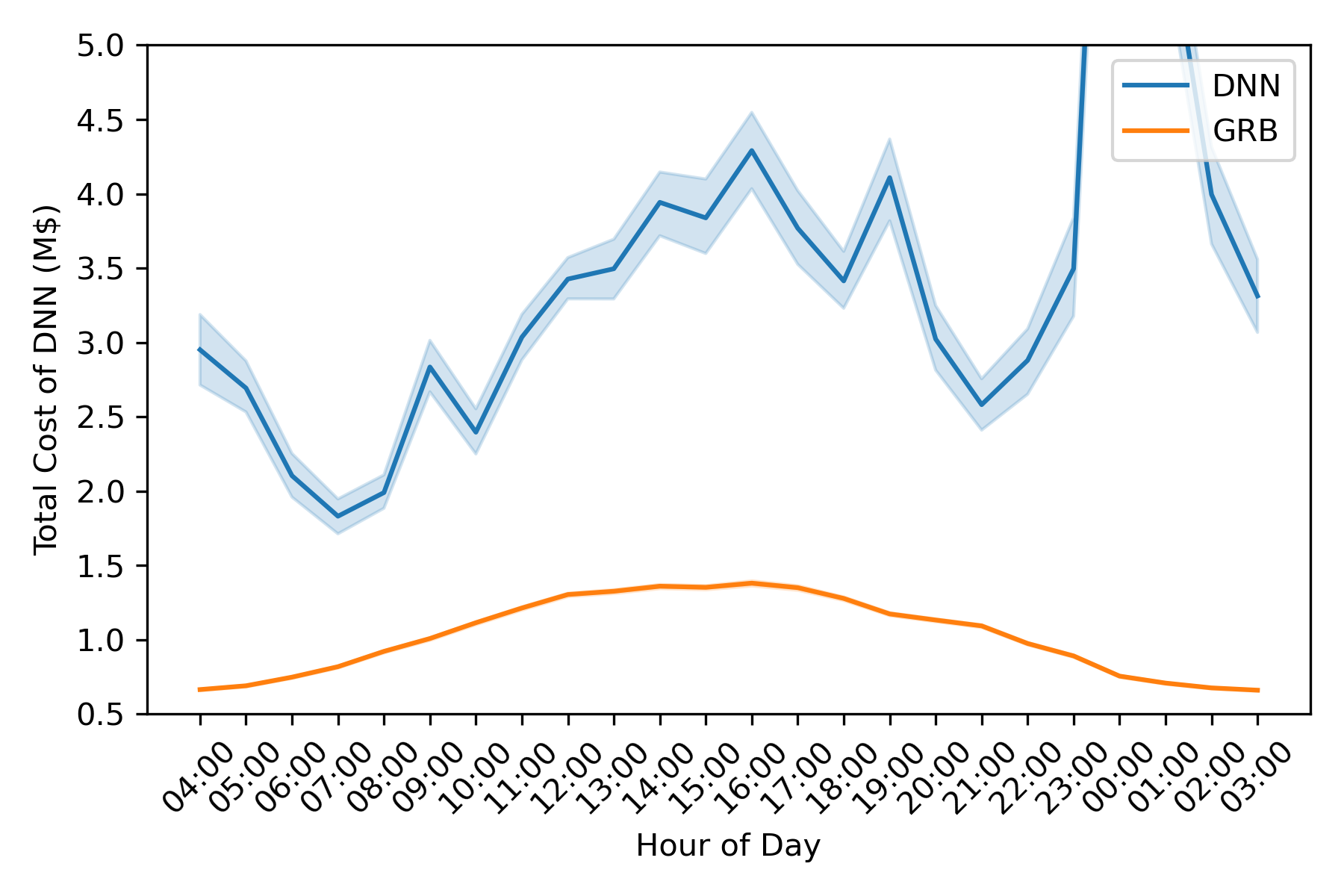}
        \hfill
        \includegraphics[width=0.44\columnwidth]{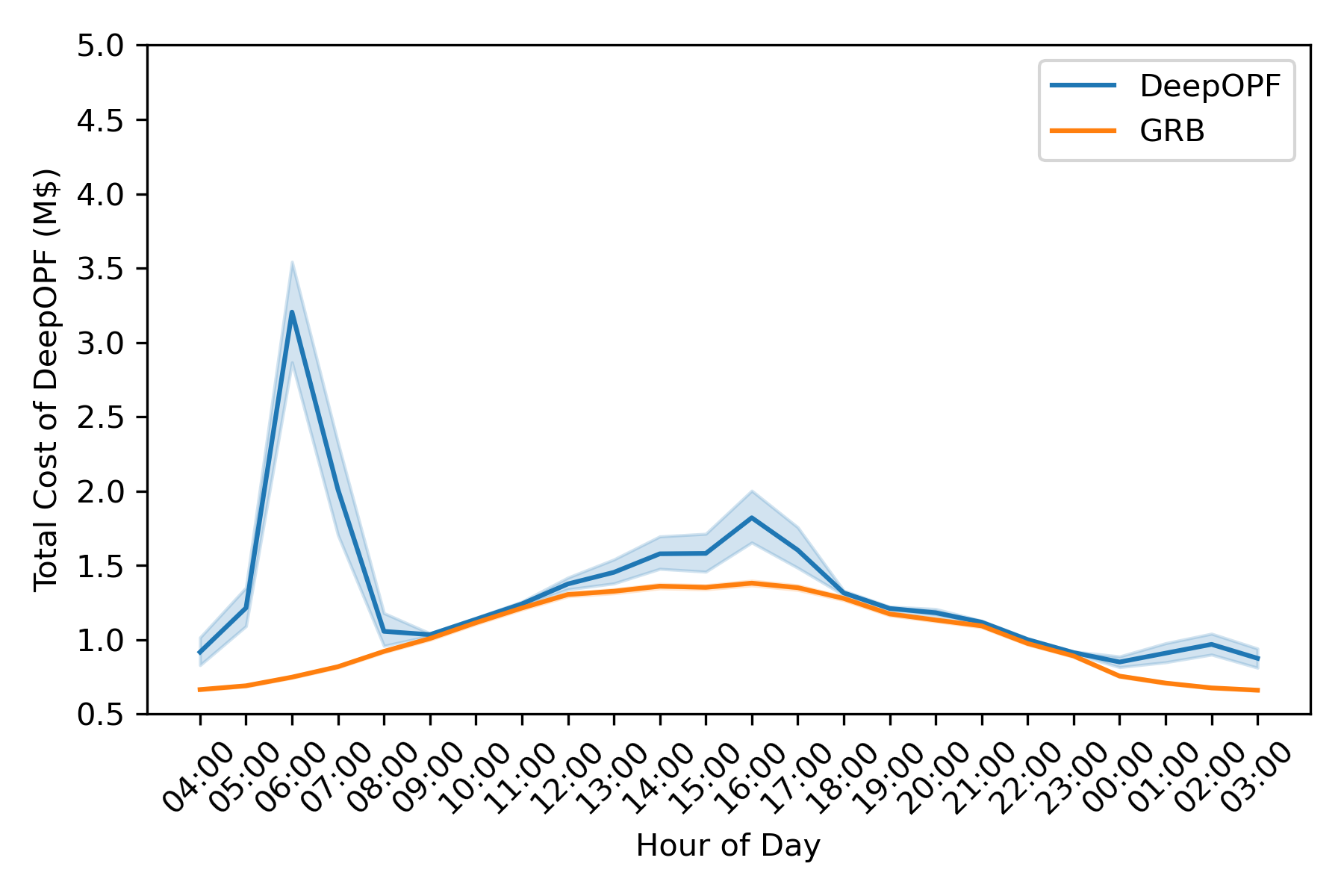}\\
        \includegraphics[width=0.44\columnwidth]{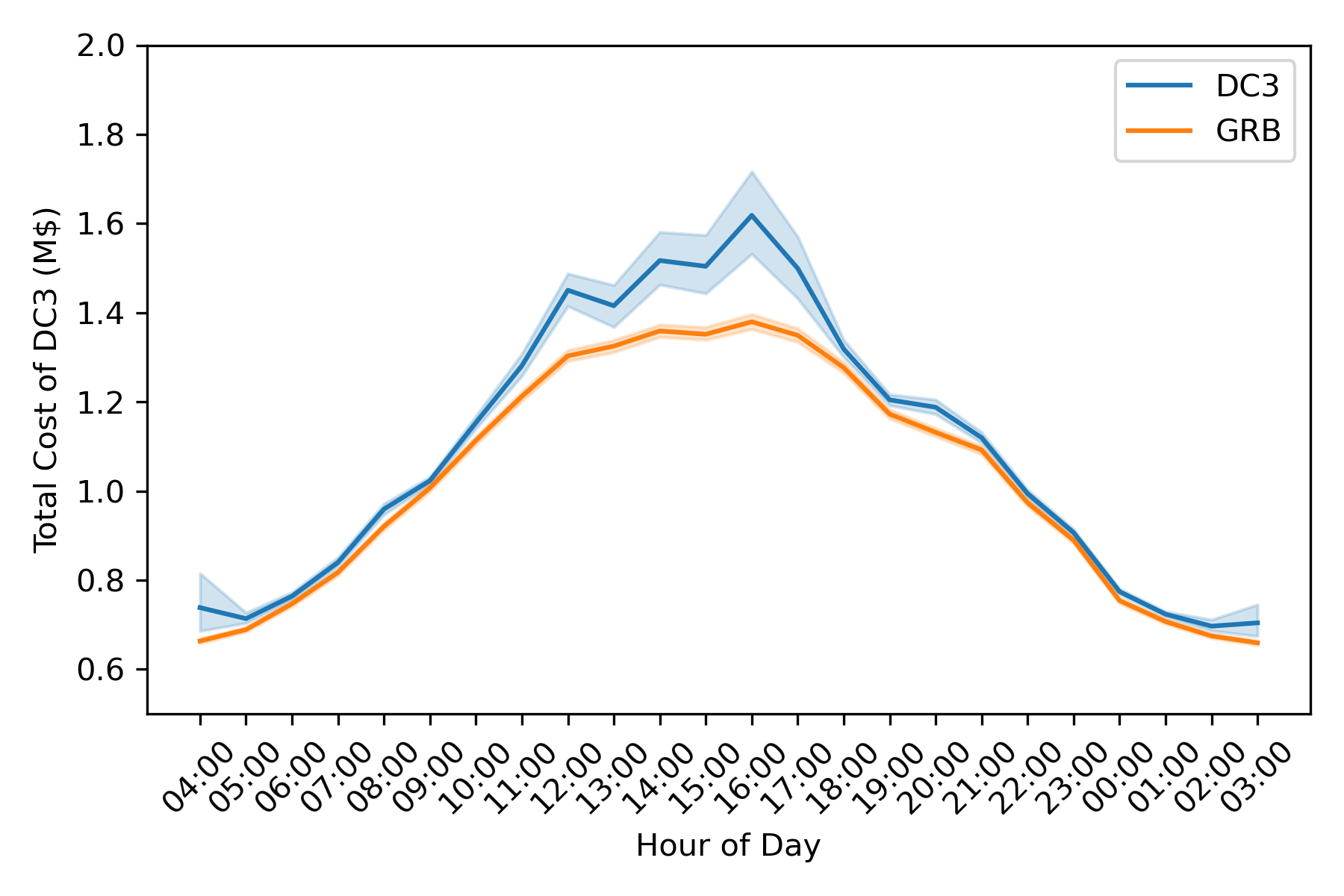}
        \hfill
        \includegraphics[width=0.44\columnwidth]{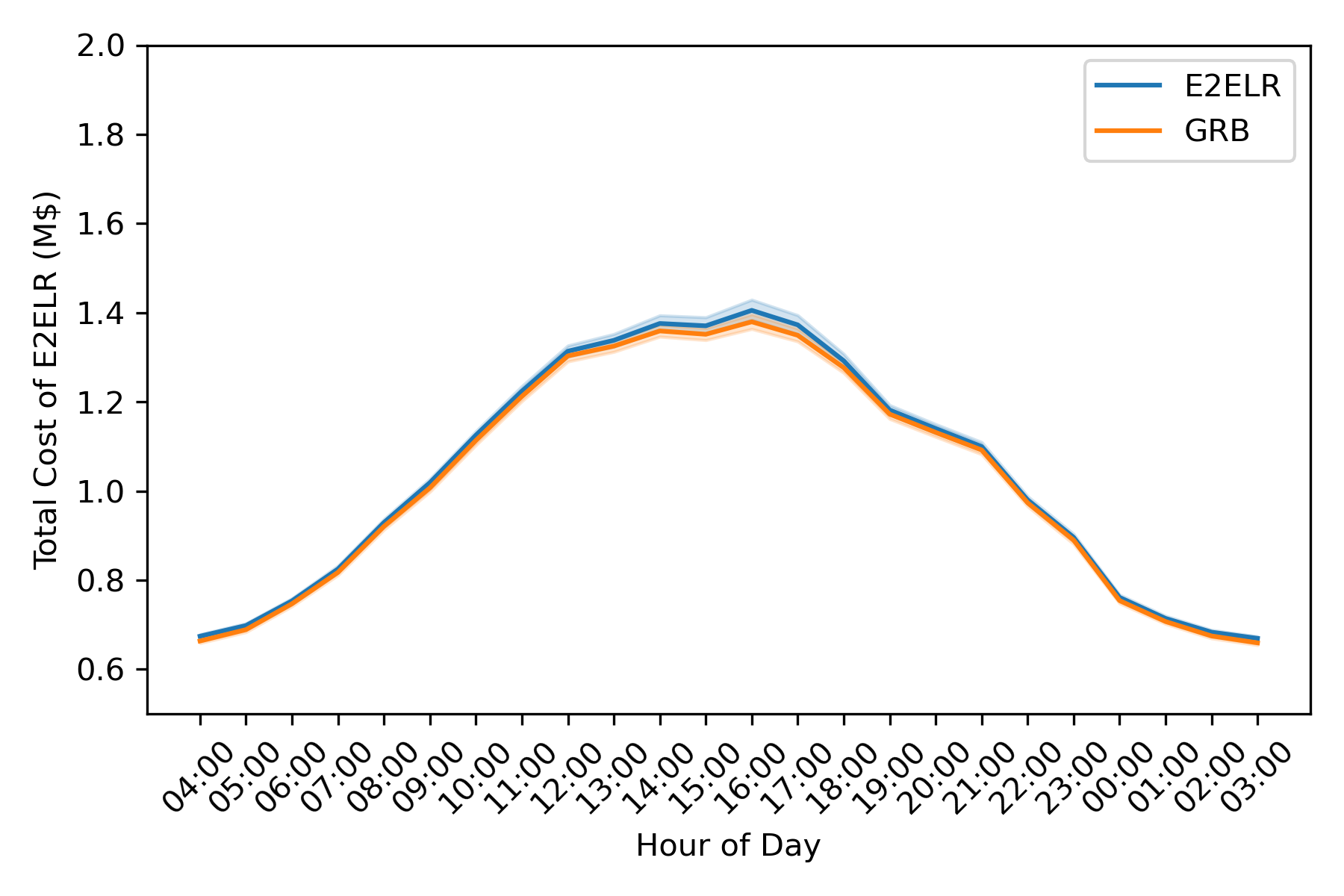}
        \caption{Hourly distributions of total cost on \pegase: comparison between ground truth (GRB) and proxies. Top: RA-DNN (left) and RA-DeepOPF (right); Bottom: RA-DC3 (left) and RA-E2ELR (right). All values in M\$.}
        \label{fig:risk:1354_pegase:overall_cost}
    \end{figure}

In addition to operational risk, grid operators rely on accurate cost
estimates to operate the grid as economically as possible.  Therefore,
Figure \ref{fig:risk:1354_pegase:overall_cost} depicts, for each proxy
architecture, the hourly distribution of total costs for the \pegase{}
system.  Here, total costs include generation costs, thermal violation
costs, and power balance violation costs.  Each figure displays the
2.5\%, mean and 97.5\% quantiles of ground truth cost distribution (in
orange), and the cost distribution obtained by the proxy (in blue).

    First, throughout the entire day, RA-DNN significantly misestimates total costs.
    Similarly, RA-DeepOPF over-estimates costs costs, especially in the early morning (4am--8am), afternoon (2pm--5pm) and at night (12am--3am).
    Both of these are caused by incorrect predictions of system imbalance and thermal violations (see Figure \ref{fig:risk:1354_pegase}).
    While RA-DC3 is more accurate, it still over-estimates total costs between 10am and 7pm, also mainly because of errors in power balance violations.
    The RA-E2ELR architecture achieves the best results, and almost perfectly matches the ground truth cost distribution.
    
    % \begin{figure}[!t]
    %      \centering
    %     \includegraphics[width=0.44\columnwidth]{img/300_ieee_compiled/300_ieee_CVAR_90_cost.png}
    %     \hfill
    %     \includegraphics[width=0.44\columnwidth]{img/1354_pegase_compiled/1354_pegase_CVAR_90_cost.png}\\
    %     \includegraphics[width=0.44\columnwidth]{img/300_ieee_compiled/300_ieee_CVAR_0_cost.png}
    %     \hfill
    %     \includegraphics[width=0.44\columnwidth]{img/1354_pegase_compiled/1354_pegase_CVAR_0_cost.png}
    %     \caption{CVAR (top) and average (bottom) of the total cost on 300 IEEE (left) and 1354 Pegase system (right). DNN is omitted due to its poor performance.}
    %     \label{fig:1354_pegase:overall_cost}
    % \end{figure}
    
\begin{figure}[!t]
        \centering
        \includegraphics[width=0.44\columnwidth]{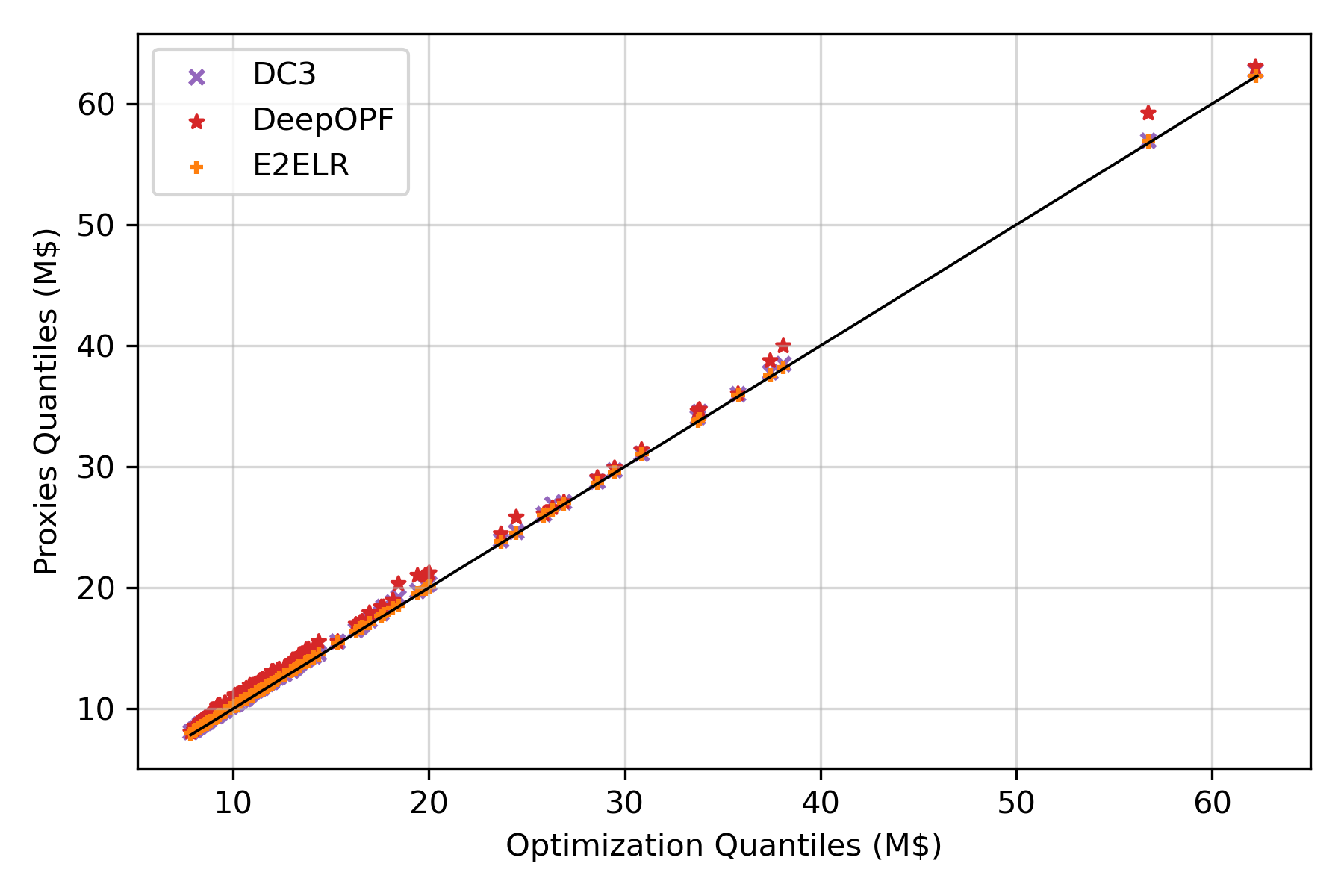}
        \hfill
        \includegraphics[width=0.44\columnwidth]{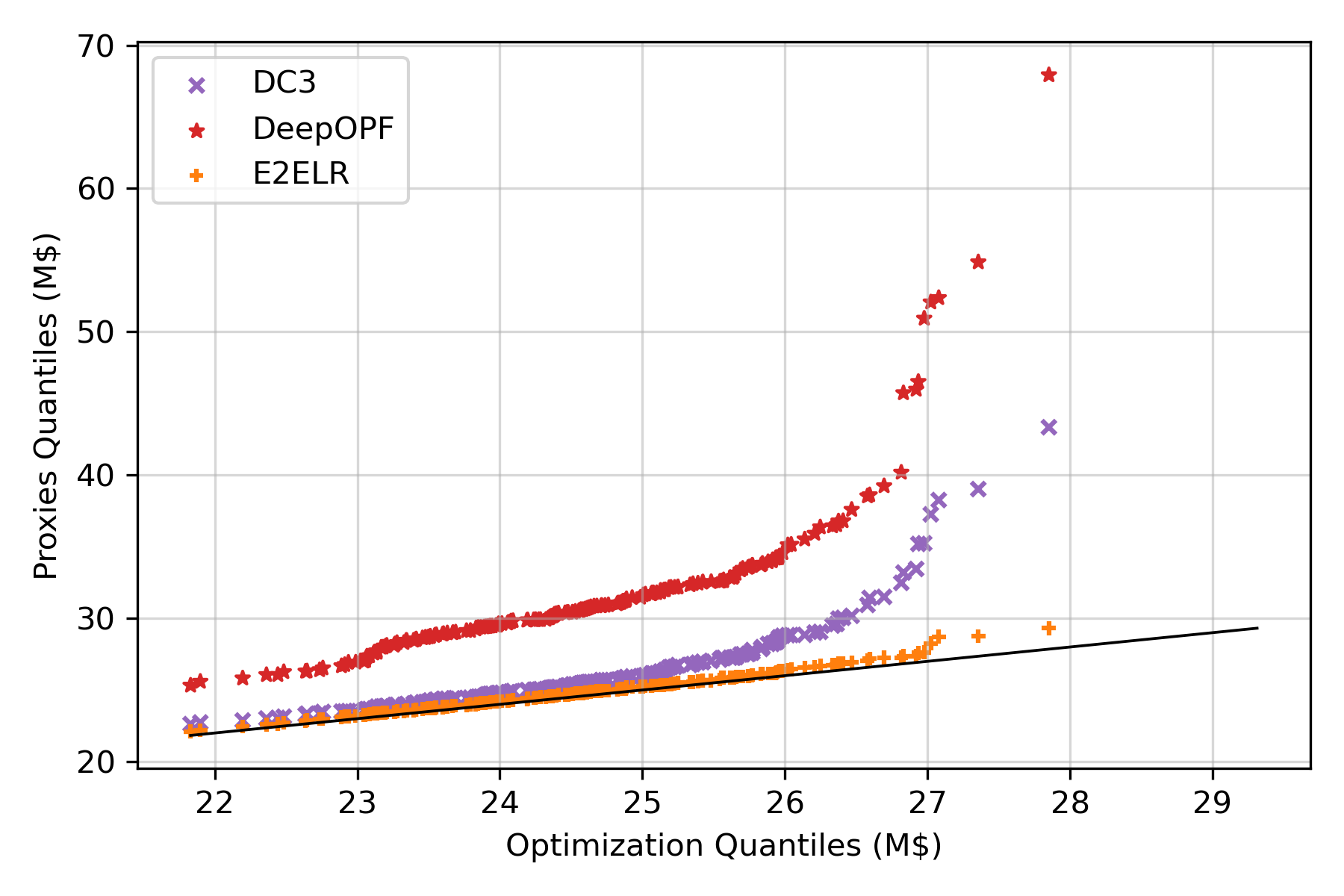}
        \caption{QQ plot of the total cost on 300 IEEE (left) and 1354 Pegase system (right). RA-DNN is omitted due to its poor performance.}
        \vspace{-1em}
        \label{fig:risk:QQ}
    \end{figure}

To complement the above analysis, Figure \ref{fig:risk:QQ} presents
quantile-quantile (QQ) plots that compare, for both systems, the
distribution of total costs (over the entire day) obtained by the
ground truth optimization-based and proxy-based simulations. RA-DNN is
omitted from the figure due to its poor performance.  For \ieeeSmall,
all proxies accurately capture the distribution of total costs.
However, for \pegase, RA-DeepOPF and RA-DC3 significantly deviate from the
ground truth distribution: both architectures significantly
over-estimate total costs for costlier (and riskier) scenarios.  This
is mostly due to incorrect predictions of system imbalances, which
incur high penalty costs.  In contrast, RA-E2ELR matches the ground truth
distribution almost perfectly.
    % These observations are consistent with the previous analyses.

\subsection{Component-level Risk Analysis}

\begin{figure}[!t]
        \centering
        \subfloat[The IEEE 300-bus system]{
            \includegraphics[width=0.44\columnwidth]{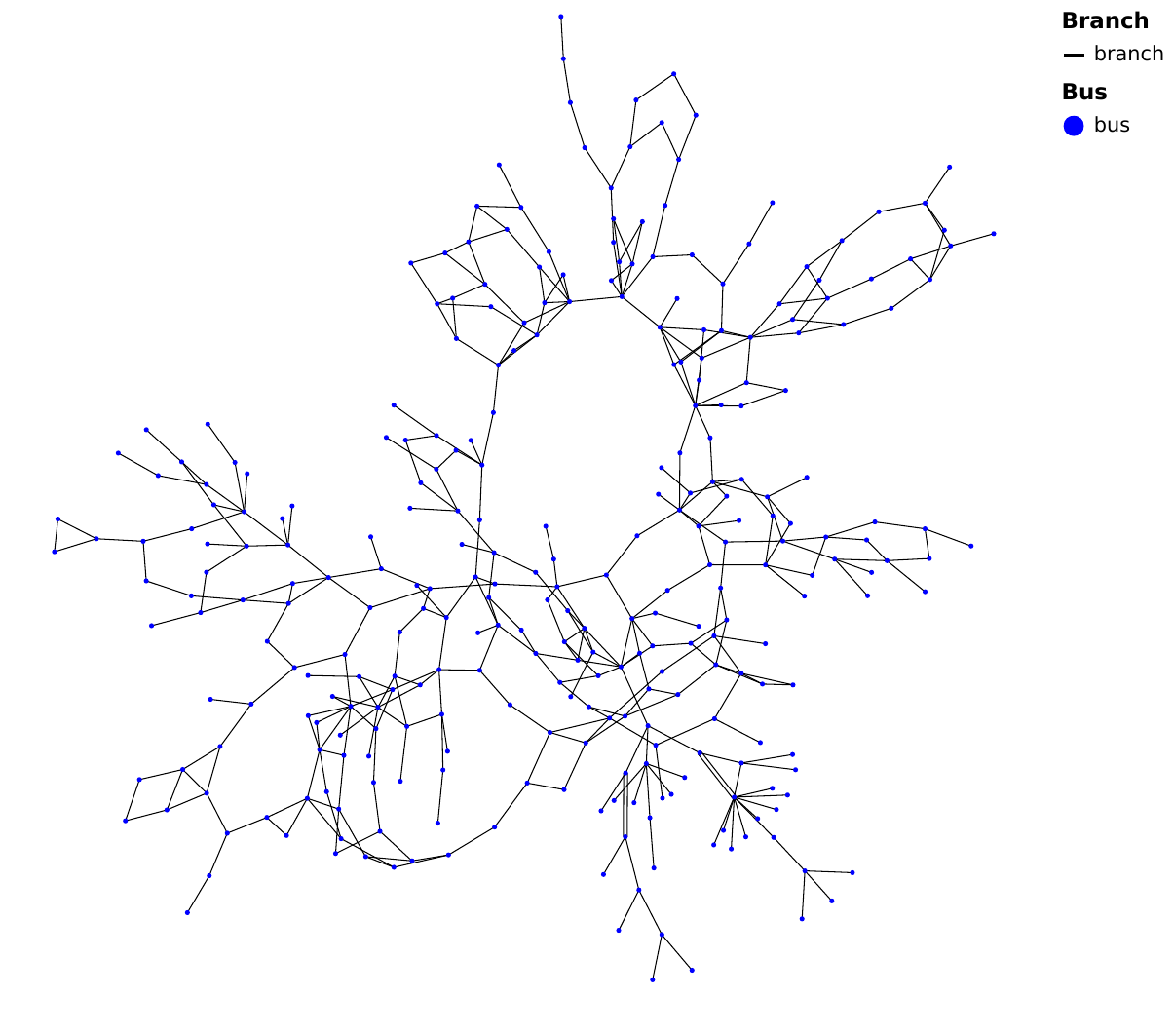}
            \label{fig:risk:branch_level_violations:grid}
        }
        \hfill
        \subfloat[Branch-level probability of thermal violations (ground truth)]{
            \includegraphics[width=0.44\columnwidth]{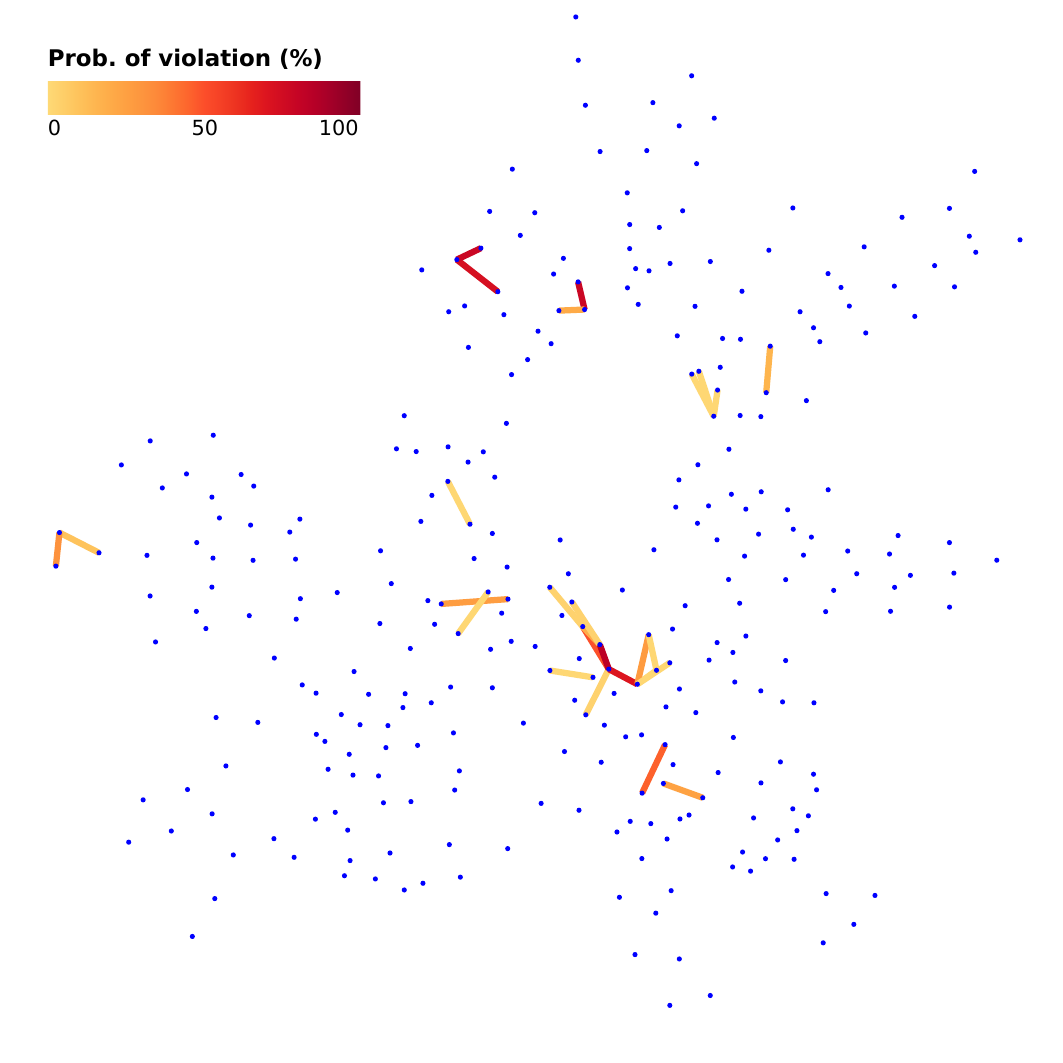}
            \label{fig:risk:branch_level_violations:GRB}
        }\\
        \subfloat[Branch-level probability of thermal violation (RA-DC3)]{
            \includegraphics[width=0.44\columnwidth]{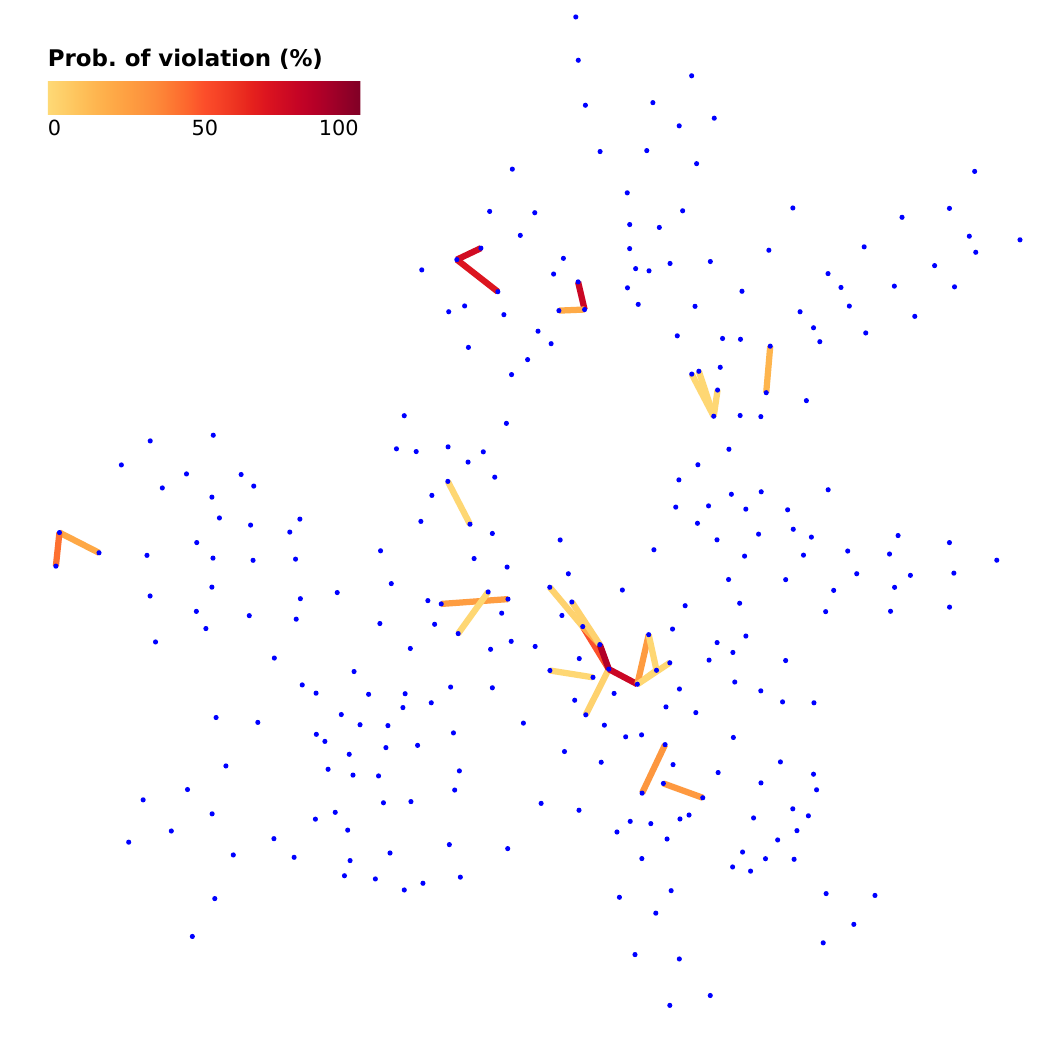}
            \label{fig:risk:branch_level_violations:DC3}
        }
        \hfill
        \subfloat[Branch-level probability of thermal violation (RA-E2ELR)]{
            \includegraphics[width=0.44\columnwidth]{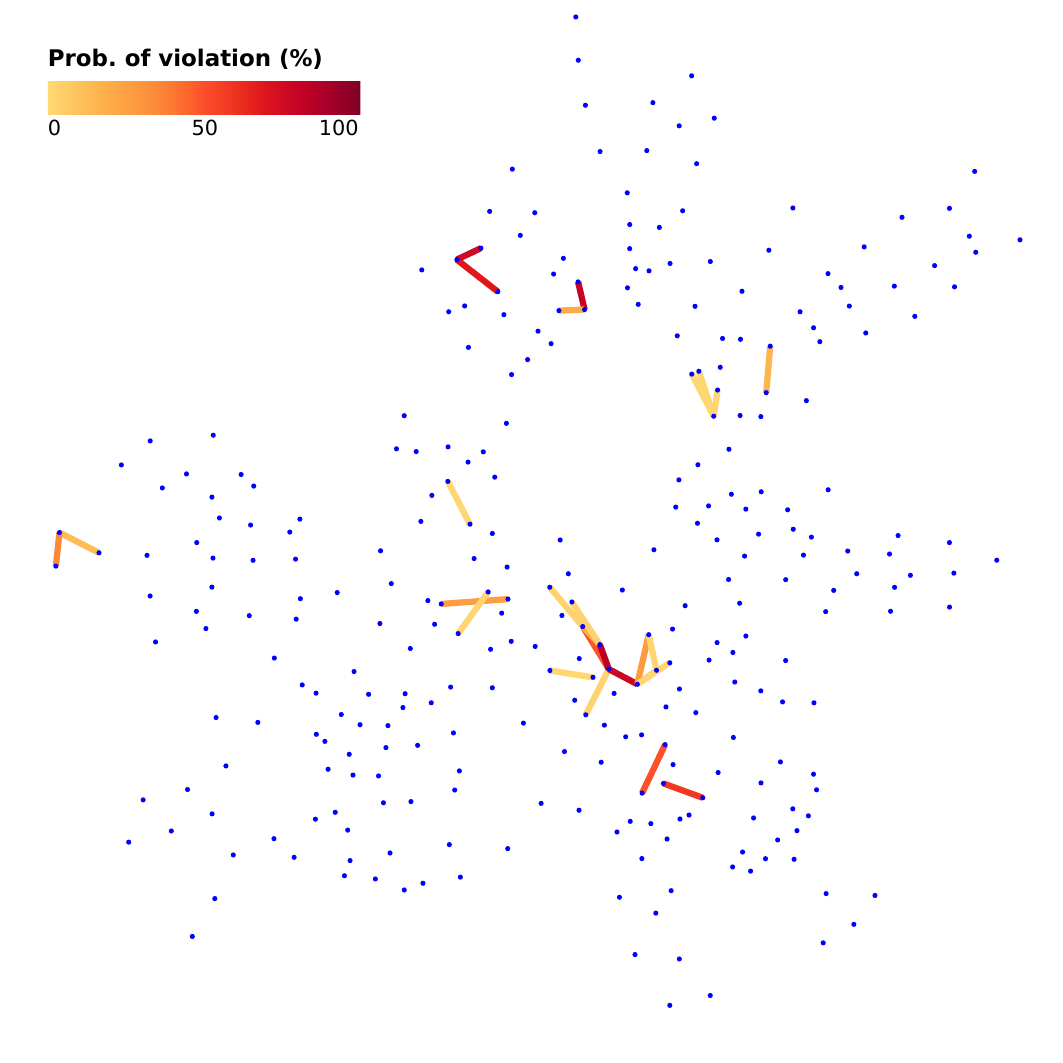}
            \label{fig:risk:branch_level_violations:E2ELR}
        }
        \caption{Branch-level probability of thermal violations for the \ieeeSmall{} system at 12pm.}
        \label{fig:risk:branch_level_violations}
        \vspace{-0.5em}
    \end{figure}

One of the advantages of the proposed proxy-based risk assessment
frameworks is its granularity, i.e., it can produce risk profiles down
to individual components.  To demonstrate this capability, Figure
\ref{fig:risk:branch_level_violations} compares ground truth and
proxy-based probability of thermal violations on each individual
branch in \ieeeSmall{} at 12pm, which is the hour with the most
violations. Figure \ref{fig:risk:branch_level_violations:grid}
displays the original grid, Figure
\ref{fig:risk:branch_level_violations:GRB} displays the ground truth
result, and Figures \ref{fig:risk:branch_level_violations:DC3} and
\ref{fig:risk:branch_level_violations:E2ELR} display the proxy-based
results obtained with RA-DC3 and RA-E2ELR, respectively.  The proxies
correctly identify all branches with non-zero probability of thermal
violations.

\subsection{Computing Times}

Table \ref{tab:exp:simulation_time} compares computing times for
proxy-based and optimization-based (GRB) risk assessment methods.
Computing times for GRB scale roughly linearly with the size of the
system, taking an average 0.5 and 2.2s per simulation per CPU core.
Note that congested scenarios typically lead to increased computing
times.  In contrast, the simulation time of the proxies is consistent
across different systems, typically taking less than 100 milliseconds
to run 100 simulations.  Specifically, RA-E2ELR achieves 30.3x speed up
than GRB, even under the assumption that hundreds of ED instances can
be solved in parallel.  Among the proxies, RA-DC3 is the slowest
because every inference requires 200 gradient steps.

    \begin{table}[!t]
        \centering
        \caption{Risk simulation computing times}
        \label{tab:exp:simulation_time}
        \begin{tabular}{lrrrrr}
            \toprule
            Systems & \multicolumn{1}{l}{E2ELR} & \multicolumn{1}{l}{DC3$^{\dagger}$} & \multicolumn{1}{l}{DeepOPF} & \multicolumn{1}{l}{DNN} & GRB$^{*}$ \\
            \midrule
             \ieeeSmall{} & 71.5 & 192.0 & 77.6 & 50.2 &  488.0 \\
             \pegase{} & 72.0 & 211.4 & 59.0 & 49.6 & 2188.0 \\
            \bottomrule
        \end{tabular}\\
        \footnotesize{All times in ms. ML inference times are per batch of 100 scenarios on 1 GPU. $^{*}$Average time per scenario using 1 CPU. $^{\dagger}$with 200 gradient steps. }
        \vspace{-1em}
    \end{table}

\section{Conclusion}

    The paper has presented the first risk-assessment framework for power systems, RA-E2ELR, that uses end-to-end feasible optimization proxies.
    The proposed framework leverages the speed of optimization proxies and the feasibility guarantees of the E2ELR architecture.
    Numerical experiments demonstrated that RA-E2ELR yields high-quality risk estimates, and a 30x speedup over optimization-based simulations.
    This contrasts with previously-proposed architectures like DeepOPF and DC3, which offer no significant speed advantage over E2ELR but whose lack of feasibility guarantees results in over-predicting risks and costs.

\bibliographystyle{IEEEtran}
\bibliography{refs}

% \end{thebibliography}

\end{document}